\documentclass[a4paper,11pt]{amsart}
\newcommand{\emailadd}{2} 
\usepackage{t1enc}
\def\frak{\mathfrak}
\def\Bbb{\mathbb}
\def\Cal{\mathcal}

\newcommand{\newc}{\newcommand}

\let\ccdot\cdot
\def\cdot{\hbox to 2.5pt{\hss$\ccdot$\hss}}

\newcommand{\om}{\omega}
\renewcommand{\phi}{\varphi}
\newcommand{\ph}{\varphi}

\newcommand{\si}{\sigma}

\newc{\aR}{\mbox{\boldmath{$ R$}}}
\newc{\aS}{\mbox{\boldmath{$ S$}}}
\newc{\aT}{\mbox{\boldmath{$ T$}}}
\newc{\aW}{\mbox{\boldmath{$ W$}}}

\newcommand{\aB}{\mbox{\boldmath{$ B$}}}

\newc{\aDeR}{\mbox{\boldmath{$ U$}}_B{}^P{}_C{}^Q}
\newc{\al}{\mbox{\boldmath$ \Delta$}}

\newcommand{\afl}{\mbox{$
\begin{picture}(9,8)(1.6,0.15)
\put(1,0.2){\mbox{$ \al \hspace{-7.8pt} /$}}
\end{picture}$}}                                           

\newcommand{\fb}  {\mbox{$                                      
\begin{picture}(9,8)(1.6,0.15)
\put(1,0.2){\mbox{$ \Box \hspace{-7.8pt} /$}}
\end{picture}$}} 

\newc{\nda}{\mbox{\boldmath$ \nabla$}}
\newc{\BigO}{\mathit{O}}

\newcommand{\aRic}{\mbox{\boldmath{$ \Ric$}}}

\newc{\aK}{\mbox{\boldmath{$ K$}}}
\newc{\aL}{\mbox{\boldmath{$ L$}}}

\newtheorem{theorem}{Theorem}[section]
\newtheorem{lemma}[theorem]{Lemma}
\newtheorem{proposition}[theorem]{Proposition}

\newcommand{\cH}{{\Cal H}}

\newcommand{\ce}{{\Cal E}}
\newcommand{\cp}{{\Cal P}}
\newcommand{\cq}{{\Cal Q}}

\newcommand{\ct}{{\Cal T}}

\newcommand{\bR}{{\Bbb R}}
\newcommand{\bT}{{\Bbb T}}

\newcommand{\bY}{{\Bbb Y}}

\usepackage{amssymb}
\usepackage{amscd}

\newcommand{\nd}{\nabla}

\newcommand{\Rho}{\mathit{P}}
\newcommand{\Up}{\Upsilon}

\newcommand{\End}{\operatorname{End}}

\newcommand{\Ric}{\operatorname{Ric}}

\newcommand{\tstar}{\mbox{\Large $ \star$}}
\newcommand{\mstar}{\mbox{\large $ \star$}}
\newcommand{\sbg}{\mbox{\scriptsize\boldmath{$g$}}}
\newcommand{\cN}{\Cal N}

\let\d=\delta
\let\hash=\sharp  
\let\#=\sharp        

\newcommand{\bA}{{\Bbb A}}

\newcommand{\bB}{{\Bbb B}}
\newcommand{\bK}{{\Bbb K}}
\newcommand{\bM}{{\Bbb M}}
\newcommand{\bU}{{\Bbb U}}

\newcommand{\upl}                         
{\mbox{$                                  
\begin{picture}(12.7,8)(-.5,-1)
\put(2.15,2.8){\line(1,0){8}}
\put(6.15,-1.2){\line(0,1){8}}
\put(6.2,2.8){\oval(8,8)[t]}
\end{picture}$}}

\def\endrk{\hbox{$|\!\!|\!\!|\!\!|\!\!|\!\!|\!\!|$}}

\newcommand{\miniX}{\mbox{\boldmath{$\scriptstyle{X}$}}}

\newcommand{\LX}{{\Cal L}_{\!\miniX}}
\newcommand{\NX}{{\mbox{\boldmath$ \nabla$}}_{\!\miniX}}

\newcommand{\aM}{\tilde{M}}
\newcommand{\ala}{\al_\alpha}
\newcommand{\act}{\mbox{\boldmath{$\ct$}}}

\def\LOT{\textit{lots}}

\newcommand{\cR}{{\Cal R}}

\newcommand{\cW}{{\Cal W}}

\newcommand{\IT}[1]{{\rm(}{\it{#1}}{\rm)}}

\newcommand{\Pa}{{\Bbb I}}

\newcommand{\nn}[1]{(\ref{#1})}

\newcommand{\bD}{\mbox{{$\Bbb D$}}}
\newcommand{\D}{\mbox{\boldmath{$ D$}}}
\newcommand{\X}{\mbox{\boldmath{$ X$}}}
\newcommand{\sX}{\mbox{\scriptsize\boldmath{$X$}}}        
\newcommand{\h}{\mbox{\boldmath{$ h$}}}
\newcommand{\bg}{\mbox{\boldmath{$ g$}}}

\newcommand{\cce}{\tilde{\ce}}                          

\newcommand{\tU}{\tilde{U}}

\newcommand{\Weyl}{C}
\newcommand{\V}{\mathit{P}}                   
\newcommand{\J}{\mathit{J}}
\newc{\CY}{A}
\newc{\Bach}{B}
\newc{\cB}{\mathcal{O}}
\newc{\SymTen}{\mbox{\sf T}}

\newcommand{\rpl}                         
{\mbox{$
\begin{picture}(12.7,8)(-.5,-1)
\put(2.3,2.7){\line(1,0){8}}
\put(6.3,-1.3){\line(0,1){8}}
\put(6.2,2.8){\oval(8,8)[r]}
\end{picture}$}}

\newcommand{\lpl}                         
{\mbox{$
\begin{picture}(12.7,8)(-.5,-1)
\put(2.3,2.7){\line(1,0){8}}
\put(6.3,-1.3){\line(0,1){8}}
\put(6.2,2.8){\oval(8,8)[l]}
\end{picture}$}}

\usepackage{ifthen}

\newc{\tensor}[1]{#1}
\newc{\Mvariable}[1]{\mbox{#1}}
\newc{\down}[1]{{}_{#1}}
\newc{\up}[1]{{}^{#1}}

\def\sideremark#1{\ifvmode\leavevmode\fi\vadjust{\vbox to0pt{\vss
 \hbox to 0pt{\hskip\hsize\hskip1em
 \vbox{\hsize3cm\tiny\raggedright\pretolerance10000
 \noindent #1\hfill}\hss}\vbox to8pt{\vfil}\vss}}}%
                        
\begin{document}

\begin{abstract} We construct here a conformally invariant
differential operator on algebraic Weyl tensors that gives special
curved analogues of certain operators related to the deformation
complex and that, upon application to the Weyl curvature, yields the
(Fefferman-Graham) ambient obstruction tensor. This new definition of
the obstruction tensor leads to simple direct proofs that the
obstruction tensor is divergence-free and vanishes identically for
conformally Einstein metrics. Our main constructions are based on the
ambient metric of Fefferman-Graham and its relation to the conformal
tractor connection. We prove that the obstruction tensor is an
obstruction to finding an ambient metric with curvature harmonic for a
certain (ambient) form Laplacian. This leads to a new ambient formula
for the obstruction in terms of a power of this form Laplacian acting
on the ambient curvature. This result leads us to construct Laplacian
type operators that generalise the conformal Laplacians of
Graham-Jenne-Mason-Sparling. We give an algorithm for calculating
explicit formulae for these operators, and this is applied to give
formulae for the obstruction tensor in dimensions 6 and 8.  As
background to these issues, we give an explicit construction of the
deformation complex in dimensions $n\geq 4$, construct two related
(detour) complexes, and establish essential properties of the
operators in these.
%
%
%
\ifthenelse{\emailadd=1}
{ 
\vspace{3mm}
\\
\newlength{\eaddwid}
\settowidth{\eaddwid}{Authors' e-mail addresses:\ \ }
\parbox[t]{\eaddwid}{Authors' e-mail addresses:}
\parbox[t]{55mm}{
\texttt{gover@math.auckland.ac.nz}\\
\texttt{lawrence.peterson@und.nodak.edu}
}
}  
{} 
%
%
\end{abstract}

\title{The ambient obstruction tensor and the conformal deformation
complex%
}

\author{A. Rod Gover and Lawrence J. Peterson}
\date{}
\address{Department of Mathematics\\
  The University of Auckland\\
  Private Bag 92019\\
  Auckland 1\\
  New Zealand} \email{gover@math.auckland.ac.nz}
\address{Department of Mathematics\\
  The University of North Dakota\\
  Grand Forks, ND 58202-8376\\
  USA} \email{lawrence.peterson@und.nodak.edu}

\maketitle

\pagestyle{myheadings}
\markboth{A.R. GOVER AND L.J. PETERSON}{THE OBSTRUCTION TENSOR AND
DEFORMATION COMPLEX}

\thanks{ARG gratefully acknowledges support from the Royal Society of 
New Zealand via Marsden Grant no.\ 02-UOA-108,
and from the New Zealand Institute of Mathematics and its Applications
for support via a Maclaurin Fellowship.}
\section{Introduction}

The Bach tensor \cite{Bach} has long been considered an
important natural invariant in 4-dimensional Riemannian and
pseudo-Riemannian geometry and continues to play an interesting role.
See \cite{Anderson,TianViac}, for example.  It is conformally invariant,
vanishes for metrics that are conformal to Einstein metrics, and
arises as the total metric variation of the action $\int |C|^2$, where
$C$ denotes the Weyl curvature.  {}From the latter and the conformal
invariance of the Weyl curvature, it follows that it is a symmetric
trace-free 2-tensor which involves 4 derivatives of the metric. An
explicit formula for the Bach tensor in terms of the Weyl curvature
$C$, the Ricci tensor, and the Levi-Civita connection is very simple:
\begin{equation}\label{one}
B_{ab}= \nd^c\nd^d C_{acbd}+\frac{1}{2}\Ric^{cd} C_{acbd}.
\end{equation}

In higher even dimensions $n$, an analogue of the Bach tensor was
discovered by Fefferman and Graham \cite{FGast}; it arose as an
obstruction to their ambient metric construction. This {\em
  Fefferman-Graham obstruction tensor}, which we denote $\cB_{ab}$ (or
sometimes $\cB^{n}_{ab}$), shares many of the properties of the Bach
tensor.  It is a trace-free symmetric 2-tensor that vanishes for
conformally Einstein metrics. 
The obstruction tensor has the form
$\Delta^{n/2-2} \nd^c\nd^d C_{acbd} +\LOT$. Here ``$\LOT$'' indicates
lower order terms. There is strong evidence that the obstruction tensor
will be as important in each even dimension as the Bach tensor is in
dimension 4. Very recently Graham and Hirachi \cite{GrH} have shown
that $\cB_{ab}$ is the total metric variation of $\int Q$, where $Q$ is
Branson's Q-curvature \cite{tomsharp,tbbo91}. This generalises the
situation in dimension 4, since in that case $\int Q$ and $\int |C|^2$
agree up to a multiple. There is a direct link between the
obstruction tensor and the non-existence of certain operators on
conformal manifolds which also generalises the 4-dimensional setting
\cite{GoH} and further indicates the critical role of the obstruction
tensor.

Despite this progress, the obstruction tensor has remained somewhat
mysterious, partly due to the lack of a general formula. In the next
section we explain that there is a fundamental difference between the
Bach tensor in dimension 4 and the obstruction tensor in even
dimensions 6 and greater. The idea is as follows. {}From the Bianchi
identities, the expression \nn{one} for the Bach tensor can be written
as $\nd^{(c}\nd^{d)} C_{acbd}+\frac{1}{2}\Ric^{cd} C_{acbd}$, where
the parentheses indicate symmetrisation over the index pair $cd$. 
The differential operator $\nd^{(c}\nd^{d)} +\frac{1}{2}\Ric^{cd} $ is
a {\em conformally invariant operator} which acts on the bundle of
``algebraic Weyl tensors'' (i.e. the bundle whose sections are 4-tensor
fields with the same conformal weight and algebraic symmetries as the
Weyl curvature) and takes values in a (density weighted) irreducible
tensor bundle. One might hope that a similar result would hold in
higher dimensions. This is not the case. In
Proposition~\ref{nonexist}, we establish that in dimensions $n\geq 6$,
the obstruction tensor cannot arise in this manner from a conformally
invariant operator that acts between irreducible tensor bundles.  This
is an easy consequence of representation theory results of
Boe-Collingwood \cite{BC} which give a classification of conformally
invariant operators on the sphere.  (See \cite{EastSlo} and references
therein.)  One focus of this article is to describe the correct
generalisation of the described construction of the Bach tensor. This
is Theorem \ref{obstrnthm}, which is one of the main results.

In the conformally flat setting, the conformally invariant operator
defined in the previous paragraph is the formal adjoint of an operator
in the so-called (conformal) deformation complex. This is a complex of
conformally invariant differential operators arising in connection
with infinitesimal deformations of a conformal structure based at a
conformally flat metric. The linearisation of the obstruction tensor,
which we denote $\sf B$, is an operator in a class of conformally
invariant operators acting between bundles in the complex. These
``long operators'' are predicted by the Boe-Collingwood
classification.  In Proposition~\ref{factoring}, we show that the
linearised obstruction operator and another long operator, that we
denote $\sf L$, factor through operators from the complex. For
example, we obtain that ${\sf B}={\sf G}{\sf C}$, where $\sf C$ is the
linearised Weyl curvature operator and $\sf G$ is a gauge companion
operator for $\sf L$. That is, $\sf L$ and $\sf G$ have the same
domain space (algebraic Weyl tensors), the system $({\sf L},{\sf G})$
gives a conformally invariant equation, and in Riemannian signature
this system is elliptic.  Theorem \ref{obstrnthm} gives a curved
analogue of this picture. The theorem describes a conformally
invariant differential operator $\bB$ which, on general conformal
manifolds of even dimension, acts on algebraic Weyl tensors and takes
values in a reducible bundle.  In dimensions $n\geq 6$, composing this
with projection to a quotient gives a conformally invariant operator
$L$ which takes algebraic Weyl tensors to weighted algebraic Weyl
tensors; $L$ generalises $\sf L$ to conformally curved manifolds. This
operator annihilates the Weyl curvature $C$, and $\bB(C)$ is the
obstruction tensor. An application of these results is given in
Proposition \ref{harmonics}, which relates (in the conformally flat
setting) the conformally invariant null space of the system $({\sf
L},{\sf G})$ to the cohomology of the deformation complex.

For a conformal structure of dimension $n$, the ambient metric is an
associated, suitably homogeneous, and Ricci-flat metric on an
$(n+2)$-manifold. In \cite{FGast}, $\cB_{ab}$ arose as an obstruction,
in even dimensions, to the existence of a formal power series solution
for this ambient metric.  In Section~\ref{invariants}, we show that
the obstruction tensor may equivalently be viewed as a formal
obstruction to having the ambient curvature harmonic for a certain
ambient form-Laplacian $\afl$.  This leads to a new proof that the
obstruction tensor is an obstruction to the ambient metric (see \IT{v}
of Theorem \ref{main1}) and a very simple ambient formula for the
obstruction.  Let $\aR$ denote the curvature of the ambient metric.
Then $\afl^{n/2-2}\aR$ is a disguised form of the obstruction. This is
also established in Theorem \ref{main1} and in the same place used to
give a new proof that the obstruction is divergence-free, i.e.\ that
$\nd^a\cB_{ab}=0$. (An alternative proof of this last result is given
in \cite{GrH}, and it also follows from the variational
characterisation given in \cite{GrH}.  See \cite{BrSrni}.)

Interpretation of these results on the underlying conformal manifold
is achieved via tractor bundles. The standard tractor bundle is a
vector bundle with a conformally invariant connection that we may view
as arising as an induced structure from the Cartan bundle with its
normal conformal Cartan connection.  On the other hand, this rank
$n+2$ vector bundle also arises in a simple way from the tangent
bundle of the ambient manifold.  Using this observation, we construct
(Theorem~\ref{powerslap} and Proposition~\ref{powerslaptan}) families
of conformally invariant operators with leading term a power of the
Laplacian; these act between arbitrary tractor bundles of an
appropriate density weight and generalise the GJMS operators of
\cite{GJMS}.  In Theorem~\ref{main2}, we show that the obstruction
tensor is obtained by applying one of these operators, namely
$\fb_{n/2-2}$, which has the form $\Delta^{n/2-2} +\LOT$, to the
tractor field $W$ that corresponds to $\aR$. Thus the problem of
finding formulae for the obstruction tensor is reduced to
understanding the special case $\fb_{n/2-2}$ of the generalised
GJMS-type operators $\fb_k$.

There is a 1-1 correspondence between between Einstein metrics and a
class of parallel standard tractors \cite{Gau,GoNur}. 
This, with the tractor formula for the
obstruction $\fb_{n/2-2}W$, forms the basis of the proof of Theorem
\ref{ceinv}, which shows that the obstruction vanishes for conformally
Einstein metrics.

Theorem~\ref{powerslap} constructs a very general class of Laplace
type conformal operators.  The inductive steps leading to
Theorem~\ref{powerslap} yield a simple and effective algorithm for
calculating explicit formulae for the conformal Laplacian operators of
that theorem. Hence by Theorem~\ref{main2}, they give an algorithm for
calculating explicit formulae for the obstruction. This algorithm is
efficient in the sense that it does not entail constructing the
ambient manifold but uses just its existence; the algorithm recovers
only those invariants of the ambient metric that actually turn up in
the final formula for the operator.  In Section~\ref{examples},
explicit tractor formulae for conformal Laplacian operators are given.
See expressions \nn{6Lap} and \nn{L9May4a}. These are then applied to
the $W$-tractor to give formulae for the obstruction in dimensions 6
and 8. Tractor formulae are given in \nn{BSix} and \nn{obsdim8}, and
formulae in terms of the Levi-Civita connection and its curvature are
given in \nn{B6} and in Figure~\ref{B8}.

The next section establishes the basic background and notation before
constructing the conformal deformation complex and introducing some
related operators. It is a pleasure to thank Tom Branson and Robin
Graham for helpful discussions.

\section{Relationship to the conformal deformation complex} \label{def}

We first sketch here notation and background for conformal structures.
Further details may be found in \cite{CapGoamb,GoPet} or
\cite{BrGodeRham}.  We mainly follow the notational conventions of the
last of these.  Let $M$ be a smooth manifold of dimension $n\geq
3$. To simplify our discussions we assume $M$ is orientable.  Recall
that a {\em conformal structure\/} on $M$ is a smooth ray subbundle
$\cq\subset S^2T^*M$ whose fibre over $x$ consists of conformally
related metrics at the point $x$.
The principal bundle $\pi:\cq\to M$ has structure group $\Bbb R_+$,
and so each representation ${\Bbb R}_+ \ni x\mapsto x^{-w/2}\in {\rm
End}(\Bbb R)$ induces a natural line bundle on $ (M,[g])$ that we term
the conformal density bundle $E[w]$. We shall write $ \ce[w]$ for the
space of sections of this bundle. Here and throughout the article,
sections, tensors, and functions are always smooth.  When no confusion
is likely to arise, we will use the same notation for a bundle and its
section space.

We write $\bg$ for the {\em conformal metric}, that is the
tautological section of $S^2T^*M\otimes E[2]$ determined by the
conformal structure. This will be used to identify $TM$ with
$T^*M[2]$.
For many calculations we will use abstract indices in an obvious way.  
Given a choice of
metric $ g$ from the conformal class, we write $ \nabla$ for the
corresponding Levi-Civita connection. With these conventions the
Laplacian $ \Delta$ is given by $\Delta=\bg^{ab}\nd_a\nd_b=
\nd^b\nd_b\,$.  
Note $E[w]$ is trivialised by a choice of metric $g$ from the
conformal class, and we write $\nd$ for the connection corresponding
to this trivialisation.
It follows immediately that (the coupled) $ \nd_a$  preserves the 
conformal metric.  

The curvature $R_{ab}{}^c{}_d$ of the Levi-Civita connection is known
as the Riemannian curvature and is defined by $
[\nd_a,\nd_b]v^c=R_{ab}{}^c{}_dv^d.  $ (Here and below,
$[\cdot,\cdot]$ indicates the usual commutator bracket.)  The
Riemannian curvature can be decomposed into the totally trace-free
Weyl curvature $C_{abcd}$ and a remaining part described by the
symmetric {\em Schouten tensor} $\Rho_{ab}$, according to $
R_{abcd}=C_{abcd}+2\bg_{c[a}\Rho_{b]d}+2\bg_{d[b}\Rho_{a]c}, $ where
$[\cdots]$ indicates the antisymmetrisation over the enclosed indices.
The Schouten tensor is a trace modification of the Ricci tensor
$\Ric_{ab}$ and vice versa: $\Ric_{ab}=(n-2)\Rho_{ab}+\J\bg_{ab}$,
where we write $ \J$ for the trace $ \V_a{}^{a}$ of $ \V$.
Under a {\em conformal transformation} we replace a choice of metric $
g$ by the metric $ \hat{g}=e^{2\om} g$, where $\omega$ is a smooth
function. Explicit formulae for the corresponding transformation of
the Levi-Civita connection and its curvatures are given in e.g.\ 
\cite{BEGo,GoPet}. We recall that in particular the Weyl curvature is
conformally invariant $\widehat{C}_{abcd}=C_{abcd}$.

A notion that we will use later is that of total order of a tensor.  A
tensor $T^{a\cdots b}{}_{c\cdots d}$ of weight $w$ and with $k$
contravariant indices and $\ell$ covariant indices will be said to be
of {\em total order} $\ell-k-w$. For example, the Weyl curvature, the
Schouten tensor, and the scalar curvature all have total order 2. The
conformal metric $\bg_{ab}$ has total order zero, and so the total
order of any tensor is unchanged by the raising and lowering of
indices using the conformal metric.

We will be interested in certain natural differential operators.  We
say that a differential operator $P$ is a {\em natural differential
operator} if it can be written as a universal polynomial in
covariant derivatives with coefficients depending polynomially on the
 metric, its inverse, the curvature tensor, and its covariant
derivatives.  The coefficients of natural operators are called {\em
natural tensors}. In the case that they are scalar they are often also
called {\em Riemannian invariants}. Note that for any tensor $T$ with
total order $t$, $\nd T$ has total order $t+1$. It follows immediately that 
for any natural differential operator $P$ that has $T$ in its domain,
the total order of $P T$ is at least $t$. 
 We say $P$ is a {\em conformally
invariant differential operator} if it is well-defined on conformal
structures (i.e.\ is independent of a choice of conformal scale).

We will use $E^k$ as a convenient alternative notation for
$\wedge^kT^*M$. The tensor product of $E^k\otimes E^\ell$, $ \ell\leq
n/2$, $k\leq \lceil n/2 \rceil$, decomposes into irreducibles. We
denote the highest weight component by $E^{k,\ell}$.  (Here ``weight''
does not refer to conformal weight, but rather the weight of the
inducing O$(n)$-representation.) We realise the tensors of
$E^{k,\ell}$ as trace-free covariant $(k+\ell)$-tensors $T_{a_1\cdots
a_kb_1\cdots b_\ell}$ which are skew on the indices $a_1\cdots a_k$
and also on the set $b_1\cdots b_\ell$. Skewing over more than $k$
indices annihilates $T$, as does symmetrising over any 3 indices.
Then we write, for example, $E^{k,\ell}[w]$ as a shorthand for the
tensor product $E^{k,\ell}\otimes E[w]$. The space of sections of each
of these bundles is indicated by replacing $E$ with $\ce$.
These sections are the {\em algebraic Weyl tensors} as discussed in the
introduction, that is, tensors $u_{abcd}$ with the same symmetries and
weight as the Weyl curvature.  In particular, the Weyl curvature itself is a
section in $\ce^{2,2}[2]$.  We will also often use the notation
$E_{k,\ell}[w]$ as a shorthand for $E^{k,\ell}[w+2k+2\ell-n]$. This
notation is suggested by the duality between $\ce^{k,\ell}[w]$ and
$\ce_{k,\ell}[-w]$; for $\phi \in \ce^{k,\ell}[w] $ and $\psi \in
\ce_{k,\ell}[-w] $, with one of these compactly supported, there is
the natural conformally invariant global pairing
\[
\phi,\psi \mapsto \langle \phi,\psi \rangle :=\int_M \phi\cdot \psi\,
d\mu_{\sbg},
\]
where $\ph\cdot \psi\in\ce[-n]$
denotes a complete contraction between $\phi $ and $\psi$.

Since the Weyl curvature is conformally invariant, it follows easily
that the linearisation (at a conformally flat metric) of the
non-linear operator $g\mapsto C^{g}\in \ce^{2,2}[2]$ (with $C^g$ the
Weyl curvature of the metric $g$) is a conformally invariant operator
${\sf C}:\ce^{1,1}[2]\to \ce^{2,2}[2]$. The formal adjoint of a
conformally invariant operator is again conformally invariant.  In
particular, the formal adjoint of $\sf C$ is conformally invariant:
$$
{\sf C}^*: \ce_{2,2}[-2] \to \ce_{1,1}[-2].
$$
Now observe that in dimension 4 we have
$\ce^{2,2}[2]=\ce_{2,2}[-2]$, and so ${\sf C}^*$ acts on the space
$\ce^{2,2}[2] $, i.e.\ the algebraic Weyl tensors.  It is given
explicitly (up to a multiple) by $U_{abcd}\mapsto
(\nd^{(a}\nd^{c)}+\Rho^{ac})U_{abcd}$.  It is straightforward to
verify directly, using the transformation formulae from \cite{GoPet},
that this is also conformally invariant in the general curved case (or
alternatively this is immediate from \nn{DIform}), and this operator
applied to the Weyl curvature gives the Bach tensor.

On conformally flat structures of dimension at least 4, the null
space of $\sf C$ locally agrees with the range of the conformal
Killing operator ${\sf K}:\ce^{1}[2]\to \ce^{1,1}[2]$ given by
$v_a\mapsto \nd_{(a}v_{b)_0}$ (where $(\cdots)_0$ indicates the
symmetric trace-free part). These operators give the initial sequence
of the {\em conformal deformation complex}. On oriented structures of
dimension 4 this complex is simply
$$
\ce^1[2] \stackrel{\sf K}{\to}\ce^{1,1}[2]\stackrel{\sf C}{\to} 
\ce^{2,2}[2] 
\stackrel{{\sf C}^*\star }{\to} \ce_{1,1}[-2] \stackrel{{\sf K}^*}{\to} \ce_1[-2],
$$
where $\tstar$ is the (conformal) Hodge star operator. Recall that
in even dimensions this gives an isomorphism on the space of middle
forms $\tstar:\ce^{n/2}\to \ce^{n/2}$, and so it also gives an isomorphism 
$\tstar:\ce^{n/2,2}[2]\to \ce^{n/2,2}[2]$.

The situation is more complicated in higher dimensions. In the 
deformation complex, the operator $\sf C$ is  followed by
the Weyl-Bianchi operator ${\sf Bi}:\ce^{2,2}[2]\to \ce^{3,2}[2]$,
given (in a conformal scale) by
\begin{equation}\label{WeylBi}
U_{abcd}\mapsto (n-3)\nd_{[a}U_{bc]de}-
g_{d[a}\nd_{|s|}U_{bc]}{}^s{}_e+ g_{e[a}\nd_{|s|}U_{bc]}{}^s{}_d. 
\end{equation}
Here the vertical bars $|\cdot|$ indicate that the enclosed indices
are omitted from the skew symmetrisation process.  (Note that an easy
consequence of its symmetries is that the operator \nn{WeylBi} is
trivial in dimension 4.)  On oriented structures the formal adjoints
of these operators conclude the complex, and so we have the picture
$$
\begin{picture}(350,20)(13,-4)
\put(10,0)
{$\cdot \stackrel{\sf K}{\to}\ce^{1,1}[2]\stackrel{\sf C}{\to} 
\ce^{2,2}[2] \stackrel{\sf Bi}{\to}\ce^{3,2}[2] \to 
\cdots \to \ce_{3,2}[-2] \stackrel{\bar{\sf Bi}}{\to} \ce_{2,2}[-2] 
\stackrel{{\sf C}^*}{\to} \ce_{1,1}[-2] \stackrel{{\sf K}^*}{\to } \cdot
$}                                                         
\end{picture}
$$ Here we have omitted the initial and terminal section spaces
($\ce^1[2]$ and $\ce_1[-2]$ respectively), since they are outside the
main focus of our discussions. In dimensions other than 6, $\bar{\sf
Bi}$ is ${\sf Bi}^*$. In dimension 6 it means the composition ${\sf
Bi}^*\tstar$.  The Hodge star is also implicitly used in interpreting
the diagram in dimension 5. In this case it gives isomorphisms
$\tstar:E^{2,2}[2]\to E_{3,2}[-2]$ and $\tstar: E^{3,2}[2]\to
E_{2,2}[-2]$, and under these ${\sf Bi}$ is identified, modulo a sign,
with ${\sf Bi}^*$. In the dimensions $n\geq 5$, ${\sf C}^*$ is given
by the same formula as in dimension 4, viz.\ $U_{abcd}\mapsto
(\nd^{(a}\nd^{c)}+\Rho^{ac})U_{abcd}$.  In even dimensions $n\geq 8$,
the centre of the pattern consists, in an obvious way, of operators
${\sf Bi}_{(k)}:\ce^{k,2}[2]\to \ce^{k+1,2}[2]$ for $k=3,\cdots
n/2-1$, their formal adjoints ${\sf Bi}_{(k)}^*:\ce_{k+1,2}[-2]\to
\ce_{k,2}[-2]$ for $k=3,\cdots n/2-2$, and ${\sf
Bi}_{(n/2-1)}^*\tstar:\ce^{n/2,2}[2]\to \ce_{n/2-1,2}[-2]$.  The
operators ${\sf Bi}_{(k)}$ generalise \nn{WeylBi}, which can be viewed
(up to a constant multiple) as the ``$k=2$ case''. For $U\in
\ce^{k,2}$, an explicit formula is $({\sf Bi}_{(k)}U)_{a_0a_1\cdots
a_kb_1b_2} ={\rm Proj}(\nd_{a_0}U_{a_1\cdots a_kb_1b_2})$, where Proj
is the bundle morphism which executes the projection into
$\ce^{k+1,2}[2]$. In odd dimensions $n\geq 7$, we have the operators
${\sf Bi}_{(k)}$ for $k=3,\cdots \lfloor n/2-1 \rfloor$, their formal
adjoints for $k=3,\cdots \lfloor n/2-2 \rfloor$.  (The operator ${\sf
Bi}_{(\lfloor n/2-1 \rfloor)}$ is formally self-adjoint).

In each dimension, the operators of the deformation complex are all
conformally invariant, and the complex is locally exact and extends to
give a resolution (on the sheaves of germs of smooth sections) of the
sheaf of conformal Killing fields.  This is a particular generalised
Bernstein-Gelfand-Gelfand (gBGG) resolution. These resolutions are
well understood and classified through the dual theory of generalised
Verma modules, and the explicit construction of the complex above is
an immediate consequence of the (local) uniqueness of the operators in
the relevant gBGG resolution, along with explicit verification of the
conformal invariance and non-triviality of the operators
mentioned. See \cite{GG} for an alternative construction of the complex via a
theory of overdetermined systems of partial differential equations
based around Spencer cohomology.

According to the results of \cite{BC}, in even dimensions the
operators of the deformation complex are not the only conformally
invariant operators between the bundles involved. There are also
``long operators'' $\ce^{k,\ell}[2]\to \ce_{k,\ell}[-2]$, and an
additional pair of operators about the centre of the pattern.  We
obtain the operator diagram
$$
\begin{picture}(350,50)(12,-32)
\put(10,0)
{$\cdot \stackrel{\sf K}{\to}\ce^{1,1}[2]\stackrel{\sf C}{\to} 
\ce^{2,2}[2] \stackrel{\sf Bi}{\to}\ce^{3,2}[2] \to 
\cdots \to \ce_{3,2}[-2] \stackrel{\bar{\sf Bi}}{\to} \ce_{2,2}[-2] 
\stackrel{{\sf C}^*}{\to} \ce_{1,1}[-2] \stackrel{{\sf K}^*}{\to } \cdot
$}                                                         
\put(36,-31){\line(1,0){284}}
\put(36,-31){\line(0,1){24}}                               
\put(320,-31){\vector(0,1){24}}
\put(44,-27){\scriptsize${\sf B}$}

\put(85,-23){\line(1,0){180}}
\put(85,-23){\line(0,1){16}}                               
\put(265,-23){\vector(0,1){16}}
\put(93,-19){\scriptsize$\sf L$}

\put(130,-14){\line(1,0){76}}
\put(130,-14){\line(0,1){7}}                               
\put(206,-14){\vector(0,1){7}}
\put(148,-10){\scriptsize${}$}

\end{picture}
$$
for dimensions 10 or greater.  The operators in this diagram are
unique (up to multiplying by a constant), and the diagram indicates by
arrows all the operators between the bundles explicitly presented.
Thus, by implication, all compositions vanish.  The same diagram
applies in dimensions 8 and 6 with minor adjustments. In dimension 8
there are two ``short'' operators with domain $\ce^{3,2}[2]$ and two
with range $\ce_{3,2}[-2]$. {}From these there is one non-trivial
composition $ \ce^{3,2}[2]\to \ce_{3,2}[-2]$.  Similiarly in dimension
6 we have $\tstar {\sf Bi}: \ce^{2,2}[2] \to\ce^{3,2}[2]$ and ${\sf
Bi}^*: \ce^{3,2}[2]\to \ce_{2,2}[-2] $, as well as the operators
indicated, and ${\sf L}={\sf Bi}^*{\sf Bi}$.  In dimension 4 the
corresponding diagram is
$$
\begin{picture}(350,50)(-78,-28)
\put(10,0){$
\cdot \stackrel{\sf K}{\to}\ce^{1,1}[2]  
\begin{picture}(20,20)(-5,0)
\put(3,10){\scriptsize$\sf C$}
\put(0,3){$\to$}
\put(0,-3){$\to$}
\put(-1,-13){\scriptsize$\mstar {\sf C}$}
\end{picture}
\ce^{2,2}[2] 
\begin{picture}(20,20)(-5,0)
\put(-1,10){\scriptsize${\sf C}^*\mstar$}
\put(0,3){$\to$}
\put(0,-3){$\to$}
\put(2,-13){\scriptsize$ {\sf C}^*$}
\end{picture}
\ce_{1,1}[-2] \stackrel{{\sf K}^*}{\to} \cdot
$}
\put(35,-23){\line(1,0){99}}
\put(35,-23){\line(0,1){16}}                               
\put(134,-23){\vector(0,1){16}}
\put(43,-19){\scriptsize$\sf B$}
\end{picture}
$$
and in this case ${\sf B}:= {\sf C}^*{\sf C}$.  Evidently on
even-dimensional conformally flat structures there are {\em detour
  complexes} (cf.\ \cite{BrGodeRham}), where one shortcuts the
deformation complex via a long operator. The examples relevant here are
\begin{equation}\label{Bdetour}
\ce^1[2]\stackrel{\sf K}{\to}\ \ce^{1,1}[2] \stackrel{{\sf B}}{\longrightarrow} \ce_{1,1}[-2]\stackrel{{\sf K}^*}{\to}  \ce_{1}[-2]
\end{equation}
and in dimensions $n\geq 6$,
$$
\ce^1[2]\stackrel{\sf K}{\to}\ \ce^{1,1}[2] \stackrel{\sf C}{\to} \ce^{2,2}[2]
\stackrel{{\sf L}}{\longrightarrow}
\ce_{2,2}[-2] \stackrel{{\sf C}^*}{\to}  \ce_{1,1}[-2]\stackrel{{\sf K}^*}{\to}  \ce_{1}[-2] .
$$
These have applications in constructing torsion quantities which
generalise Cheeger's de Rham half-torsion \cite{torsionprogress}.

According to \cite{FGast}, the obstruction tensor $\cB_{ab}$ is a
trace-free symmetric 2-tensor of weight $2-n$.  That is, it is a
section of $\ce_{1,1}[-2]=\ce^{1,1}[2-n]$.  {}From the general theory
in \cite{EastSlo}, we know that all the operators indicated explicitly
by arrows in the diagrams above admit curved analogues, that is,
generalisations to general conformal structures. (In fact, the
formulae given above for $\sf K$, ${\sf C}^*$, and $\sf Bi$ give
conformally invariant operators on general structures. We will
continue to use this notation for these operators even in the
conformally curved setting.)  {}From the diagrams, however, the
difference between dimension~4 and higher even dimensions is clear.
In dimension~4 there is a conformal operator $\ce^{2,2}[2]\to
\ce_{1,1}[-2]$ that yields the Bach tensor, as described above. In
higher dimensions the conformally invariant ${\sf C}^*$ does not have
$\ce^{2,2}[2]$ as domain. These observations establish the following
key point.
\begin{proposition}\label{nonexist}
In even dimensions $n\geq 6$, there can be no conformally invariant
differential operator $\ce^{2,2}[2]\to \ce_{1,1}[-2]$ that recovers
the obstruction tensor upon application to the Weyl curvature $C$.
\end{proposition}
\noindent If there were such an operator, then by Theorem~\ref{main1},
below, or by \cite{GrH}, it would necessarily have as highest order term
$\Delta^{n/2-2}\nd^a\nd^c U_{abcd}$.  Its linearisation would
therefore be an operator $\ce^{2,2}[2]\to \ce_{1,1}[-2]$.  But there
is no operator between these bundles in the diagram.

This brings us to the question of whether, in dimensions
$n\geq 6$, there can be any conformally invariant operator that yields
the obstruction tensor. We will see that there is, and we will
construct the operator.  To understand how this works, it is helpful
to expose some properties of the operators $\sf B$ and $\sf L$.
\begin{proposition}\label{factoring}
  The operators ${\sf B}:\ce^{1,1}[2]\to \ce_{1,1}[-2]$ and ${\sf L}:
  \ce^{2,2}[2]\to \ce_{2,2}[-2]$ are formally self-adjoint. In each even
  dimension $n\geq 6$ the following holds: there is a natural linear differential operator ${\sf H}:
  \ce^{2,2}[2]\to \ce_{2,2}[-2]$ such that $\sf B$ is given by the
  composition
$$
{\sf B}= {\sf C}^* {\sf H} {\sf C} ;
$$
there is a natural linear differential
operator ${\sf N}: \ce^{3,2}[2]\to \ce_{3,2}[-2]$ such that $\sf L$
is given by the composition
$$
{\sf L}={\sf Bi}^* {\sf N} {\sf Bi}.
$$
\end{proposition}
A proof of this is given in Section \ref{FormulaeSect}. The proof
there uses the geometric tools we develop shortly. The factorisations
described in the proposition can also be established via central
character arguments (and see also \cite{BrGodeRham}).  Note that $\sf
L$ is only defined in even dimensions $n\geq 6$.  In dimension~6,
${\sf N}$ is the identity.
Otherwise, from the classification of conformally invariant
operators on conformally flat manifolds, as discussed above, it
follows that the operators $\sf H$ and $\sf N$ are not conformally
invariant.

On conformally flat structures the operator ${\sf G}:= {\sf C}^*{\sf
  H}$ is not conformally invariant ($n \neq 4$).  It \textit{is},
however, conformally invariant on the range of the linearised Weyl
curvature, and we have ${\sf B}={\sf G}{\sf C}$. On the other hand,
${\sf L}$ annihilates the range of $\sf C$. The theorem below gives
special curved analogues of these operators.

We need some further notation.  On conformal manifolds of dimension
$n$ there is a natural reducible, but indecomposable, bundle $W_{2,2}$
that has the composition series $E_{2,2}[-2]\lpl E_{2,1}[-2] \lpl
E_{1,1}[-2] $. This means that $E_{1,1}[-2]$ is a (conformally
invariant) subbundle and that $E_{2,1}[-2]$ is a subbundle of the
quotient $W_{2,2}/E_{1,1}[-2]$.  The bundle $W_{2,2}$ (which is a
subbundle of a certain tractor bundle) is constructed explicitly in
proof of Theorem \ref{obstrnthm} in Section~\ref{FormulaeSect}, and given a choice
of metric $g$ from the conformal class, it decomposes as $[W_{2,2}]_g=
E_{2,2}[-2]\oplus E_{2,1}[-2] \oplus E_{1,1}[-2] $.  Let us write
${\rm I}^*$ and ${\rm P}$ for the respective canonical bundle maps
$W_{2,2}\to E_{2,2}[-2] $ and $W_{2,2}\to E_{2,2}[-2]\lpl E_{2,1}[-2]
$ (which are unique up to a constant multiple).
\begin{theorem}\label{obstrnthm}
  On conformal manifolds of even dimension $n \geq 6$ there is a
  natural non-trivial conformally invariant linear differential
  operator
$$
\bB:\ce^{2,2}[2]\to \cW_{2,2}= \ce_{2,2}[-2]\lpl \ce_{2,1}[-2]\lpl  \ce_{1,1}[-2]
$$
with the following properties:\\
\IT{i} The composition $({\rm I}^*\bB=:L):\ce^{2,2}[2]\to
\ce_{2,2}[-2] $ is a non-trivial conformally invariant differential
operator of order
$n-4$.   \\
\IT{ii} There is a linear differential operator $\overline{\bB}$ such
that $\overline{\bB}\bB=\Delta^{\ell}+\LOT$.  Thus on Riemannian
signature conformal structures, $\bB$ is graded
injectively elliptic.  \\
\IT{iii} For the Weyl curvature $C\in \ce^{2,2}[2]$ we have $\bB(C)\in
\ce_{1,1}[-2]$. The natural conformal invariant
$\cB_{ab}\in\ce_{1,1}[-2]$ given this way agrees with the obstruction tensor.
\end{theorem}
\noindent We prove the theorem in Section \ref{FormulaeSect}. Note that 
there is a degenerate version of the theorem for dimension 4; see
expression \nn{Wform} and the comments that follow it.

{}From the uniqueness of $\sf L$ it is clear that on conformally flat
manifolds $L$ recovers $\sf L$ (up to a constant multiple). However
$L$ is a special curved generalisation of $\sf L$, since the property
$L(C)=0$ generalises to arbitrary conformal structures the vanishing
of the composition ${\sf L}{\sf C}$. Since $L(C)$ vanishes, it follows
from the composition series for $\cW_{2,2}$ that the component of
$\bB(C)$ in $\ce_{2,1}[-2]$ is a natural conformal invariant. That
this also vanishes is also a special property of $\bB$ that, in a
sense, carries to general structures the non-existence of an operator
$\ce^{1,1}[2]\to \ce_{2,1}[-2]$. It follows that on conformally flat
structures the composition $\bB{\sf C}$ determines a non-trivial
operator $\ce^{1,1}[2]\to \ce_{1,1}[-2]$ which therefore agrees with
${\sf B}$. If, for each metric $g$ in the conformal class, we write
$G$ for the composition of $\bB$ followed by projection to the
component $\ce_{1,1}[-2]$ (we have such a projection since, recall,
$\cW_{2,2}$ completely decomposes, given a conformal scale), then, by
construction, $G$ is a curved analogue of the operator $\sf G$.  That
is, the restriction of $G$ to conformally flat structures is $\sf G$.
Note that $G$ has the special property that $G(C)=\cB$, and (as we
will see from the construction of $\cW_{2,2}$) although $G$ is not
conformally invariant, the conformal variation of $G$ under $g\mapsto
e^{2\om} g$ is only quadratic in $\om$. Since $\sf G$ also has this
sort of variation, this is optimal.

In the conformally flat case, it is easily shown that ${\rm P}\bB $ can
be re-expressed as a composition $\bU \sf L$.  (Here $\bU$ is the
operator \nn{trext}, below, except with $w$ set to $6-n$, and this
result follows from the non-existence of a non-trivial conformal
operator $\ce^{2,2}[2]\to \ce_{2,1}[-2]$.) 
It follows  
 from this and \IT{ii}  that (in even dimensions
$n\geq 6$) $({\sf L},{\sf G})$ is a right factor of a Laplacian. That
is, there are linear differential operators $\overline{\sf L}$ and
$\overline{\sf G}$ such that
$$
(\overline{\sf L}~,~\overline{\sf G} )\left(\begin{array}{c}\sf L\\
\sf G\end{array}\right)
=\Delta^\ell+\LOT .
$$
Since also ${\sf G}$ is conformally invariant on the null space of
$\sf L$, it follows that $\sf G$ is a conformal gauge companion
operator in the sense of \cite{BrGoopava}.  (See also
\cite{BrGodeRham}).  Thus in Riemannian signature, the operator pair
$({\sf L},{\sf G})$ is an elliptic system.  Since $\sf L$ has $\sf Bi$
as a right factor, the system $({\sf Bi},{\sf G})$ is also elliptic
and has a conformally invariant null space. Let us denote this by
$\cH^2_G$, and note that on compact manifolds, $\cH^2_G$ is
finite-dimensional.  This is closely related to the second cohomology
of the deformation complex. For example, from Proposition
\ref{factoring} and an easy adaption of the proof of Proposition 2.5
in \cite{BrGodeRham}, we obtain the following result, which suggests
that $\cH^2_G$ is a candidate for a space of conformal harmonics.
Here we write $H^i$, $i=1,2$, for the first and second cohomology
spaces in the deformation complex, and $H^1_B$ for the first
cohomology of the detour complex \nn{Bdetour}.
\begin{proposition}\label{harmonics}  
  On even-dimensional conformally flat manifolds of dimension $n\geq
  6$, there is an exact sequence:
$$
0\to H^1 \to H^1_B \to \cH^2_G\to H^2 .
$$
\end{proposition}
\noindent The maps are as follows: $\cH^2_G\to H^2$ is simply $\Phi\mapsto
[\Phi]$; $H^1_B \to \cH^2_G$ is the map on the quotient $\cN({\sf
  B})/\cR({\sf K})$ induced by the restriction of ${\sf C}$ to
$\cN({\sf B})$, the null space of $\sf B$; $H^1 \to H^1_B$ is
inclusion. There are further results concerning the relationship of
$H^1_B$ to $H^1$ and $\cH^2_G$ to $H^2$, but this will be taken up
elsewhere. (See also \cite{BrGoopava}.)

\section{The ambient construction and tractor calculus}\label{ambient}

In the subsequent sections we will explore the relationship between
the Fefferman-Graham ambient metric construction \cite{FGast} and
tractor calculus as derived in \cite{CapGoamb} and \cite{GoPet}. The
notation and conventions for the ambient metric closely follow
\cite{BrGodeRham}.

For $\pi:\cq\to M$ a conformal structure of signature $ (p,q)$, let us
use $\rho $ to denote the ${\Bbb R}_+$ action on $ \cq$ given by
$\rho(s) (x,g_x)=(x,s^2g_x)$.  An {\em ambient manifold\/} is a smooth
$(n+2)$-manifold $\aM$ endowed with a free $\Bbb R_+$--action $\rho$
and an $\Bbb R_+$--equivariant embedding $i:\cq\to\aM$.  We write
$\X\in\frak X(\aM)$ for the fundamental field generating the $\Bbb
R_+$--action.  That is, for $f\in C^\infty(\aM)$ and $ u\in \aM$, we
have $\X f(u)=(d/dt)f(\rho(e^t)u)|_{t=0}$.  For an ambient manifold
$\aM$, an {\em ambient metric\/} is a pseudo--Riemannian metric $\h$
of signature $(p+1,q+1)$ on $\aM$ satisfying the conditions: (i) $\Cal
L_{\sX}\h=2\h$, where $\Cal L_{\sX}$ denotes the Lie derivative by
$\X$; (ii) for $u=(x,g_x)\in \cq$ and $\xi,\eta\in T_u\cq$, we have
$\h(i_*\xi,i_*\eta)=g_x(\pi_*\xi,\pi_*\eta)$; and (iii) $\Ric(\h)=0$
up to the addition of terms vanishing to order $n/2-1$ if $ n$ is even
(or $\Ric(\h)=0$ to all orders if $n$ is odd); (iv)
$\h(\X,\cdot)=\frac{1}{2}d Q$ to all orders. 

If $M$ is locally conformally
flat, then there is a canonical solution to the ambient metric problem
to all orders. This is simply to take a flat ambient metric. This is
forced by (i--iii) in odd dimensions, but in even dimensions this
extends the solution (cf.\ comments in \cite{BrGodeRham}).  When
discussing the conformally flat case, we assume this solution.

We write $ \nda $ for the ambient Levi-Civita connection, and we use
uppercase abstract indices $A$, $B$, etc., for tensors on $ \aM$. The
ambient Riemann tensor will be denoted $\aR_{AB}{}^C{}_{D}$.  Since
$\LX\h=2\h$, it follows that
\renewcommand{\arraystretch}{1.5}
\begin{equation} \label{ndaX} 
\nda\X =\h,
\end{equation} 
and
\begin{equation}\label{XRQ} 
\X^A\aR_{ABCD}=0. 
\end{equation} 
Equalities without qualification, as here, indicate that the results
hold to all orders or identically on the ambient manifold.

\subsection{Tractor bundles}\label{tractor}

Let $ \cce(w)$ denote the space of functions on $ \aM$ which are
homogeneous of degree $ w\in {\Bbb R}$ with respect to the action $
\rho$. 
More generally, a tensor field $F$ on $ \aM$ is said to be {\em homogeneous
  of degree} $w$ if $\rho(s)^* F= s^w F$ (i.e.\ $ \LX F=w F$).
Just as sections of $ \ce[w]$ are equivalent to functions in $
\cce(w)|_\cq$, the restriction of a homogeneous
tensor field to $\cq$ has an interpretation on $M$.
Denote by $ \act$ the space of sections of $ T\aM$ which are
homogeneous of degree $ -1$ and write $\act(w)$ for sections in $\act\otimes
\cce(w)$, where the $\otimes$ here indicates a tensor product over
$\cce(0)$. {}From \cite{CapGoamb} we have the
following results: We may identify the standard tractor bundle $\bT$
with $T\aM|_\cq$ modulo a suitable ${\Bbb R}_+$-action so that
sections of $\bT$ are in one-one correspondence with sections in
$\act$. Thus we write $\ct$ for the space of sections of the standard
tractor bundle.  The filtration of $\bT$, which we traditionally
indicate by a composition series,
\begin{equation} \label{trcomp} 
\bT= E[1]\lpl E^1[1] \lpl E[-1] ,
\end{equation}
reflects the vertical subbundle of $T\cq$ and $T\cq$ as a subbundle of
$T\aM|_\cq$.  Then since the ambient metric $\h$ is homogeneous of
degree $2$, it descends to give a metric on $\bT$.  This is the usual
tractor metric.  Sections of $ \act$ may be characterised as those
sections of $ T\aM$ which are covariantly parallel along the integral
curves of $ \X$ (which on $ \cq$ are exactly the fibres of $\pi$).
The normal tractor connection agrees with the ambient connection as
follows.  For $ U\in \ct$, let $ \tU $ be the corresponding section of
$ \act|_{\cq}$. A tangent vector field $ \xi$ on $M$ has a lift to a
field $ \tilde{\xi}\in \act(1)$, on $ \Cal Q$, which is everywhere
tangent to ${\Cal Q}$.  This is unique up to adding $f\X$, where $
f\in \cce (0)$. We extend $ \tU$ and $\tilde{\xi} $ smoothly and
homogeneously to fields on $ \aM$ and form $ \nda_{\tilde{\xi}} \tU
|_\cq$; this section is independent of the extensions and independent of the
choice of $ \tilde{\xi}$ as a lift of $ \xi$ and is exactly 
the section of $ \act(0) |_\cq$ corresponding to $ \nd_\xi U $  where 
$\nd$ here indicates the tractor connection. 

When abstract indices are required, the section spaces of the tractor
bundle and its dual can also be denoted $\ct^A$ and $\ct_A$. A choice
of metric $g$ from the conformal class determines
\cite{BEGo,CapGoluminy} a canonical splitting of the composition
series \nn{trcomp}. Via this the semi-direct sums $\lpl$ in that
series get replaced by direct sums $\oplus$, and we introduce
$g$-dependent sections $Z^A{}^b\in\ct^{Ab}[-1]$ and $Y^A\in\ct^A[-1]$
that describe this decomposition of $\bT$ into the direct sum $
[\bT^A]_g= E[1]\oplus E_a[1] \oplus E[-1]$.  A section $V\in\ct$ then
corresponds to a triple $[V]_g=(\sigma,\mu,\rho)$ of sections from the
direct sum according to $V^A= Y^A\sigma+Z^{Ab}\mu_b+X^A\rho $, and in
these terms the tractor metric is given by $h(V,V)=\bg^{ab}\mu_a\mu_b
+2 \si \rho$.  Thus the tractor contractions of the projectors are 
\begin{equation}\label{innerprodt}
X^AY_A=1,\quad Z^{Ab}Z_{Aa}=\delta^b_a,
\end{equation}
and 0 for the other pairings.

If $\hat{Y}^A$ and $ \hat{Z}^A{}_b$ are the projectors for 
the metric $ \hat{g}=e^{2\om}g$, then we have
\begin{equation}\label{XYZtrans}
\textstyle
\begin{array}{rl}
\hat Z^{Ab}=Z^{Ab}+\Up^bX^A, &
\hat Y^A=Y^A-\Up_bZ^{Ab}-\frac12\Up_b\Up^bX^A.
\end{array}
\end{equation}
Here $\Up:=d\om$.  In terms of this splitting, determined by $g$, the
tractor connection is given by
\begin{equation}\label{connids}
\begin{array}{rcl}
\nd_aX_A=Z_{Aa}\,, &
\nd_aZ_{Ab}=-\V_{ab}X_A-Y_A\bg_{ab}\,, & \nd_aY_A=\V_{ab}Z_A{}^b .
\end{array}
\end{equation}

 We use the notation $\tilde{\bT}^\Phi$ to denote an arbitrary ambient
 tensor bundle (with $ \tilde{\bT}^0$ meaning the trivial bundle) and
 write $\act^{\Phi}(w)$, $w\in {\Bbb R}$, for the subspace of
 $\Gamma(\tilde{\bT}^\Phi)$ consisting of sections $S$ satisfying
 $\nda_{\sX}S=wS$; we will say such sections are homogeneous of 
{\em weight} $w$.  {}From the
 constructions above, it follows that the sections in
 $\act^{\Phi}(w)|_\cq$ are equivalent to sections of a tractor bundle
 that we denote $\bT^\Phi[w]$.  We write $\ct^\Phi[w]$ for the section
 space of the latter.

A basic example of interest is the bundle of $k$-form tractors
$\bT^k$, which is the $k^{\underline{\rm th}}$ exterior power of the
bundle of standard tractors. It is straightforward to verify that this
has a composition series which, in terms of section spaces, is given
by
\begin{equation} \label{formtractorcomp}
\ct^k=\Lambda^k\ct\cong\ce^{k-1}[k]\lpl\{\ce^k[k]\oplus
\ce^{k-2}[k-2]\}\lpl\ce^{k-1}[k-2].
\end{equation}
Also of direct relevance to our constructions below are the bundles
which we denote $\bT^{2,2}[w]$. For a given $w\in \bR$, $\bT^{2,2}[w]$
is the subbundle of $\bT^2\otimes \bT^2 \otimes E[w]$ consisting of
tractors of weight $w$ and Weyl tensor type symmetries (that is,
Riemann tensor type symmetries and also trace-free). We write
$\ct^{2,2}[w]$ for the section space of $\bT^{2,2}[w]$ and note that
(with notation as in Section \ref{def}) it has the composition series
\begin{equation}\label{wcomp}
\renewcommand{\arraystretch}{0.5}
\begin{array}{ccccccccc}
                                &&&&\ce^{2,2}[w+4] &&&&\\
           &              &  &     & \oplus &          &  &  &     \\
           &       & \ce^{2,1}[w+4] &   & \ce^2[w+2] &    &\ce^{2,1}[w+2]   &  &  \\
\ce^{1,1}[w+4]\hspace*{-2mm}
&
\hspace*{-2mm}
\lpl
\hspace*{-2.3mm}
&
\hspace*{-2mm}
\oplus
\hspace*{-2mm}
&
\hspace*{-2mm}
\lpl
\hspace*{-2.3mm}
&
\hspace*{-2mm}
\oplus
\hspace*{-2mm}
&
\hspace*{-2mm}
\lpl
\hspace*{-2.3mm}
&
\hspace*{-2mm}
\oplus
\hspace*{-2mm}
&
\hspace*{-2mm}
\lpl
\hspace*{-2.3mm}
&
\hspace*{-2mm}
\ce^{1,1}[w] .
\\
   &       & \ce^1[w+2] &      & \ce^{1,1}[w+2] &    & \ce^1[w]   &  & \\
   &            &    &    &   \oplus &    &   &       & \\
                                &&&& \ce[w]  &&&&\\
\end{array}
\renewcommand{\arraystretch}{1}
\end{equation}
A comment on punctuation is in order: here the columns represent
composition factors, decomposed into ${\frak so}(g)$-irreducibles, and
these are separated by $\lpl$'s which indicate the composition
structure. This series may be obtained by any ${\frak so}(n+2)$ to
${\frak so}(n)$ branching-rule algorithm or, alternatively, by simply
considering the possible contractions of the projectors
$X$, $Y$, and $Z$ into a typical element of $\ct^{2,2}[w]$.

\subsection{Operators and invariants via the ambient metric}\label{invariants}

An operator $P$ acting between ambient tensor bundles is said to be
{\em homogeneous} of weight $u\in {\Bbb R}$ if $[\NX,P]=u P$.
Operators homogeneous in this sense map homogeneous tensors of weight
$w$ to homogeneous tensors of weight $w+u$. On the other hand, a
differential operator $P$ is said to act {\em tangentially} along
$\cq$, as an operator on some domain space, if we have $PQ=QP'$ for
some operator $P'$ (or equivalently $[P,Q]=QP''$ for some $P''$).  Of
particular interest are linear differential operators $P$ which are
both homogeneous and also, on some homogeneous tensor space
$\act^\Phi(w)$ as domain, act {\em tangentially} along $\cq$. Each
such operator $P$ clearly determines a well-defined operator on $
\act^\Phi(w)|_\cq$, the restriction of the relevant homogeneous tensor
space to $\cq$, and hence determines an operator on the equivalent
weighted tractor bundle section space $\ct^\Phi[w]$. If the operator
$P$ is natural as an operator on the ambient manifold, then since the
ambient construction is not dependent on a choice of metric from the
conformal class, it follows that the induced operator on weighted
tractor fields is conformally invariant. The remaining issue is
whether this induced operator is natural for the underlying conformal
structure. For the operators we are interested in here, we solve this
by giving an algorithm for expressing the induced operator as a
formula in terms of known natural operators. This solves two problems,
since one of our aims is to obtain explicit formulae for the operators
concerned.

Before we construct examples of such operators, we require some further
background.  First note that from \nn{ndaX}, we have
\begin{equation}\label{bracketalX}
[\al,\X]=2\nda,  \quad {\rm where} \quad \al:=\nda^A\nda_A ,
\end{equation}
and $\nda_A Q=2\X_A$.  Both identities hold to all orders. Thus
$\nda_{\sX}Q=2Q$; $Q$ is homogeneous of weight 2.  A short computation
shows that if $ U$ is an ambient tensor field, then
\begin{equation}\label{delQ}
[\al,Q]U =2(n+2\nda_{\sX}+2) U
.
\end{equation}
It follows that for any positive integer $\ell$, if an ambient tensor
field $U$ is $O(Q^{\ell})$,
%
%
then $\al U$ and $\nda U$ are both
$O(Q^{\ell-1})$.

Now we define an operator that we denote
$\D$ (or $\D_A$ when indices are used).  Let
\begin{equation}\label{lab23June3g}
\D V:=\nda (n+2\NX-2) V-\X\al V,
\end{equation}
for any ambient tensor field $V$. It is readily verified that $\D$ is
homogeneous of weight $-1$.  By \nn{bracketalX} we also have the
equivalent formula
\begin{equation}\label{alternativeD}
        \D V=\nda (n+2\NX) V-\al\X V.
\end{equation}
Using either of these with the computations above, we obtain 
\[
\D QV=Q\D V+4Q\nda V ,
\]
and so $\D$ acts tangentially.  
For later use we note  that for any integer $\ell\geq 2$, if
$V$ is $O(Q^\ell)$, then $\D_A V$ is
$O(Q^{\ell-1})$.  

Since $\D$ acts tangentially on {\em any} ambient tensor bundle, it
follows that for every tractor bundle $\ct^\Phi$ and $w\in \bR$ we
obtain an operator
$$
D: \ct^\Phi[w]\to \ct \otimes \ct^\Phi[w-1]
$$
equivalent to $\D$ as an operator $\act^\Phi(w)|_\cq \to \act\otimes
\act^\Phi(w-1)|_\cq$.  It is straightforward to prove (see
\cite{CapGoamb,GoPet}) that $D$ is the usual tractor-D operator of
\cite{T,BEGo}. For a given choice of metric $g$ from the conformal
class and for any $V\in\ct^\Phi[w]$, $D$ is given explicitly by
\begin{equation}\label{Dform}
D^A V:=(n+2w-2)w Y^A V+ (n+2w-2)Z^{Aa}\nabla_a V -X^A\Box V,
\end{equation}
where $\Box V:= \Delta V+ w \J V$.  We note that $D$ is a natural
differential operator.  A differential operator taking values in a
tractor bundle (or acting between tractor bundles) is said to be
natural if the ${\frak so}(g)$-irreducible components of the operator
are natural.

Note  that acting on $\act^\Phi(1-n/2)$, $\D$ is simply
$-\X\al$, and correspondingly $D$ simplifies to $-X\Box$ on
$\ct^\Phi[1-n/2]$. Thus $\al$ acts tangentially on $\act^\Phi(1-n/2)$
and, as an operator on the restriction of this space to $\cq$, is
equivalent to the tractor-coupled conformal Laplacian 
\begin{equation}\label{clapn}
 \Box:
\ct^\Phi[1-n/2] \to \ct^\Phi[-1-n/2] .
\end{equation}

Many identities involving $D$ are obtained most easily by calculating
with $\D$ on $\aM$. For example, a short calculation using \nn{ndaX}
and \nn{bracketalX} shows that
\begin{equation}\label{DAboldXA} 
  \D_A\X^AV=(n+2w+2)(n+w)V-Q\al V
\end{equation}
for any $V\in\act^{\Phi}(w)$. Hence 
for any $V\in\ct^{\Phi}[w]$, we have
\begin{equation}\label{Lab13July3d}
D_A X^A V=(n+2w+2)(n+w)V.
\end{equation}

An observation key to the next section is that the ambient curvature
$\aR$ is, at low orders at least, ``harmonic'' for a certain
Laplacian. Before we construct this Laplacian we need some further
notation.  Let us write $\# $ ({\em hash}) for the natural tensorial
action of sections $A$ of $\End(T\aM)$ on ambient tensors. For example,
on an ambient covariant 2-tensor $T_{AB}$, we have
$$
A\hash T_{AB}=-A^{C}{}_AT_{CB}-A^{C}{}_BT_{AC}.
$$
If $A$ is skew for $\h$, then at each point, $A$ is
$\frak{so}(\h)$-valued.  The hash action thus commutes with the
raising and lowering of indices and preserves the ${\rm
SO}(\h)$-decomposition of tensors. (For example, $A\#$ maps trace-free
symmetric tensors to trace-free symmetric tensors). As a section of
the tensor square of the $\h$-skew bundle endomorphisms of $T\tilde
M$, the ambient curvature has a double hash action on ambient tensors;
we write $\aR\hash \hash T$. As a point on punctuation, it should be
noted that we will treat tensors in composite expressions as
multiplication operators.  A composition of operators $L$, $M$, and $N$
acting on $S$ denoted $LMN S$ means $L(M(N(S)))$.  For example, $\nda
\aR\hash \hash T$ has the same interpretation as $\nda (\aR\hash \hash
T)$.

{}From the Bianchi identities, we have that on any
Riemannian or pseudo-Riemannian manifold,
\begin{equation}\label{dispone}
\begin{array}{lll}
\lefteqn{4\nda_{A_1}\nda_{B_1}\aRic_{A_2 B_2} =}&&\vspace{2mm}
\\ &&
\al \aR_{A_1A_2B_1B_2} +\frac{1}{2}\aR\#\#\aR_{A_1A_2B_1B_2}
\\ && 
-\aRic_{CA_1}\aR^{C}{}_{A_2B_1B_2} +\aRic_{CB_1}\aR^{C}{}_{B_2A_1A_2}
.
\end{array}
\end{equation}
\noindent{\bf Remark:}
In \nn{dispone} we adopt the convention that sequentially labelled
indices in the subscript position (such as $A_1$ and $A_2$) are
implicitly skew-symmetrised.  This convention applies throughout this
paper unless noted otherwise. \quad \endrk

Let us define a Laplacian operator $\afl$ by the formula
\[
\afl := \al +\frac{1}{2}\aR\#\# .
\]
Then from (\ref{dispone}) and the conditions on $\Ric(\h)$ for the
ambient metric, we have
\begin{equation}\label{semihar}
  4\nda_{A_1}\nda_{B_1}\aRic_{A_2 B_2}  =\afl\aR_{A_1A_2B_1B_2}
  +O(Q^{n/2-1})
\end{equation}
in even dimensions.  Therefore
\begin{equation}\label{aflaR}
\afl\aR_{BCDE}= 0
\end{equation}
modulo $O(Q^{n/2-3})$ in even dimensions and to infinite
order in odd dimensions. \\
\noindent{\bf Remarks:} 1. The operator $\afl$ is a type of form-Laplacian. 
On a Riemannian or pseudo-Riemannian manifold, suppose $U$ is any
tensor with Riemann tensor type symmetries.  A short calculation shows
that \newc{\ad}[1]{\mbox{\boldmath$ d_{#1}^{\mbox{\tiny$\nda$}}$}}
\newc{\da}[1]{\mbox{\boldmath$ \delta_{#1}^{\mbox{\tiny$\nda$}}$}}
$$ \afl U =
-\frac{1}{2}\left(\da{1}\ad{1}+\ad{1}\da{1}+\da{2}\ad{2}+\ad{2}\da{2}
\right) U,
$$
where $\ad{i}$ is the Levi-Civita connection-coupled exterior
derivative, $\da{i}$ its formal adjoint, and the index $i$ is 1 or 2
according to whether we regard $U$ as a 2-form (with values in a
tensor bundle) on the first pair of indices or the last pair. (In
terms of the Levi-Civita connection $\nda$, we have, for example,
$({\ad{1}U})_{A_0A_1A_2B_1B_2}=3\nda_{A_0}U_{A_1A_2B_1B_2}$ and
$(\da{2}U)_{A_1A_2B_2}=-\nda^{B_1}U_{A_1A_2B_1B_2}$.)

Returning to the ambient manifold, note that from these
observations, the results concerning the degree to which the ambient 
curvature is $\afl$-harmonic are
manifest, since on the one hand $\ad{1}$ and $\ad{2}$ annihilate $\aR $
by the Bianchi identity and on the other hand $\da{1}\aR$ and $\da{2}\aR$
are $O(Q^{n/2-2})$ (or $O(Q^\infty)$ in odd dimensions) by dint of the
contracted Bianchi identity and the condition (iii) on the ambient
Ricci curvature. 

\vspace{1mm}

\noindent\label{harmony}2. Note that from \nn{dispone}, if $\aRic$
vanishes to all orders on the ambient manifold, then it is immediate
that $\afl \aR$ vanishes to all orders. Conversely, if $\afl \aR$
vanishes to all orders, then so does $ 4
\nda_{A_1}\nda_{B_1}\aRic_{B_2 A_2}
+\aRic_{CA_1}\aR^{C}{}_{A_2B_1B_2}-\aRic_{CB_1}\aR^{C}{}_{B_2A_1A_2}
$.  \ On the other hand, contracting the latter with
$\X^{A_1}\X^{B_1}$ and using \nn{ndaX} and \nn{XRQ} yields $2
\aRic_{A_2 B_2}$.  Thus on the ambient manifold, the vanishing of
$\aRic$ to all orders is equivalent to the vanishing of $\afl \aR$ to
all orders.  $\quad\endrk$

\smallskip

We may view the operator $\afl$  as the special case 
$\alpha=1/2$ of the family of ambient Laplacians 
\begin{equation}\label{aLapl}
\ala := \al +\alpha\aR\#\# , \quad \alpha\in {\Bbb R},
\end{equation}
which also includes the ambient form Laplacian at $\alpha =1$ and the
usual ambient Bochner Laplacian at $\alpha =0$. While the latter was
used in the constructions of \cite{GJMS} giving conformal operators
between densities, the generalisation to the ambient form Laplacian
proved appropriate in \cite{BrGodeRham} for the study of conformal
operators on (weighted) differential forms.  It seems likely that
others in the family will also have important roles, and so much of
the discussion in the next section allows for the possibility of any
$\alpha\in {\Bbb R}$.  Certain key identities for $\al$ are unaffected
by the addition of the $\aR\#\#$ term. In particular, since
$\X^A\aR_{ABCD}=0$ it follows that
\begin{equation}\label{abracketalX}
[\ala,\X] = [\al,\X]=2\nda .
\end{equation}
Using this, or even  more simply by noting that  $[\aR\#\#,Q]=0$, we obtain
\begin{equation}\label{abracketQ}
[\ala,Q] = [\al,Q] =2(n+2\nda_{\sX}+2) .
\end{equation}
A point of departure is $[\ala,\nda]$. 
Observe that if $V_{BC\cdots E}$ is any ambient tensor,
then by the Ricci flatness of the ambient metric,
\begin{equation}\label{L17May01c}
\begin{array}{rll}
\lefteqn{[\al, \nda_A]V_{BC\cdots E}=}&&
\vspace{1mm}
\\
&& -2\aR_A{}^P{}_B{}^{Q}\nda_PV_{QC\cdots
  E}-2\aR_A{}^P{}_C{}^{Q}\nda_PV_{BQ\cdots E} -\cdots
\vspace{1mm}
\\
&&
-2 \aR_A{}^P{}_E{}^{Q}\nda_PV_{BC\cdots  Q}
.
\end{array}
\end{equation}
This equality holds modulo $O(Q^{n/2-2})$ in even dimensions and
to infinite order in odd dimensions.  

Using the results above and
the Bianchi identities, it is straightforward to verify that if we
define the ambient homogeneous (of weight $-2$) tensor field
\begin{equation}\label{WR}
\aW_{A_1A_2B_1B_2}:= \frac{3}{n-2}\D^{A_0}\X_{A_0}\aR_{A_1A_2B_1B_2}
,
\end{equation}
then in dimensions other than 4, we have
\[
\aW|_\cq=(n-4) \aR|_\cq.
\]
Note that $\aW$ is
well-defined in all dimensions and by construction is conformally
invariant.  Thus the equivalent tractor field $W_{ABCE}$ is
conformally invariant and of weight $-2$. In dimensions other than 4,
it is immediate that this has Weyl tensor type symmetries.  (Recall that
$\aR|_\cq$ is trace-free.)  In fact, it has these symmetries in all
dimensions  and is a natural tractor field.
In a choice of conformal scale, $W_{ABCE}$ is given by
\begin{equation}\label{Wform}
\begin{array}{l}
(n-4)\left( Z_A{}^aZ_B{}^bZ_C{}^cZ_E{}^e C_{abce}
-2 Z_A{}^aZ_B{}^bX_{[C}Z_{E]}{}^e A_{eab}\right. \\ 
\left.-2 X_{[A}Z_{B]}{}^b Z_C{}^cZ_E{}^e A_{bce} \right)
+ 4 X_{[A}Z_{B]}{}^b X_{[C} Z_{E]}{}^e B_{eb},
\end{array}
\end{equation}
where $A_{abc} $ is the {\em Cotton tensor},
\begin{equation}\label{cot}
A_{abc}:=2\nabla_{[b}\Rho_{c]a} ,
\end{equation} 
 and 
\begin{equation}\label{Bachform}
B_{ab}:=\nabla^c
A_{acb}+\Rho^{dc}C_{dacb}.
\end{equation}
Note that from \nn{Wform} it follows that, in dimension 4, $B_{eb}$ is
conformally invariant. 
This is the Bach tensor: 
from the contracted Bianchi identity, we have
\begin{equation}\label{bi1}
(n-3)A_{abc}=\nabla^d C_{dabc} ,
\end{equation} 
and so in dimension 4 \nn{Bachform} agrees with \nn{one}. In other dimensions $n\geq 3$ we also 
refer to $B_{ab}$, as defined in
\nn{Bachform}, as the {\em Bach tensor}.  The
tractor field $W$ first appeared in \cite{gosrni,goadv}. The
connection to the ambient curvature was derived in \cite{CapGoamb},
where the above results are treated in detail.

\section{Conformal Laplacians and the ambient obstruction}\label{FormulaeSect}

In this section we show how one can obtain the ambient obstruction
tensor by applying a conformally invariant operator $\fb_{n/2-2}$ of
the form $\Delta^{n/2-2}+\LOT$ to the natural tractor field $W$ defined
above.  For any integer $m\geq 1$, we let
\[
\fb_m:=\Box_m^{1/2},
\]
where $\Box_m^{1/2}$ is the case $\alpha=1/2$ of the operator
$\Box_m^{\alpha}$ of Theorem~\ref{powerslap}, below.  We prove
Theorem~\ref{powerslap} in Section~\ref{lapsect}.  The inductive
nature of the proof of Theorem~\ref{powerslap} will show that one can
construct explicit tractor formulae for the operators $\Box_m^\alpha$
in terms of $X$, $D$, $W$, $h$, and $h^{-1}$.  One may thus use
Theorem~\ref{main2} together with a choice of conformal scale and the
formula for $W$ given in \nn{Wform} to construct a tractor formula for
$\cB_{ab}$.  It is then easy to expand this tractor formula to a
formula in terms of the Levi-Civita connection and its curvature.

In what follows, the phrase ``generic $n$-even case'' refers to the
case in which $n$ is even and $M$ is conformally curved.
\begin{theorem} \label{powerslap}
  For every integer $m\geq 1$ and for every $\alpha\in{\Bbb R}$, there
  exists a conformally invariant operator
  $\Box^\alpha_m:\ct^{\Phi}[m-n/2]\to \ct^{\Phi}[-m-n/2]$ having
  leading term $\Delta^m$ which is natural as follows: in odd
  dimensions and for conformally flat $M$ for all $m\geq 1$; in the
  generic $n$-even case for $1\leq m\leq n/2-2$, or if $\alpha=0$ for
  $1\leq m\leq n/2-1$, or if $\ct^{\Phi}[m-n/2]=\ct[m-n/2]$ for $1\leq
  m\leq n/2-1$, or if $\ct^{\Phi}[m-n/2]=\ct^0[m-n/2]$ for $1\leq
  m\leq n/2$.  In these cases there is a tractor formula for
  $\Box_m^\alpha$ which is given by a partial contraction polynomial
  in $\Box$, $D$, $W$, $X$, $h$, and $h^{-1}$, and
  this polynomial is linear in $U$. In the tractor formula for
  $\Box_m^\alpha U$, each free index appears either on $U$ or on a
  $W$-tractor.
\end{theorem}

We believe the operators $\fb_m$ will be important for many
problems. For our current purposes, we are primarily interested in
them when $n$ is even, $m=n/2-2$, and the domain bundle is
$\ct^{2,2}[-2]$.  In particular, we have the following result, which
is an immediate consequence of Theorem \ref{main1}, below.
\begin{theorem} \label{main2}
Let $M$ be a conformal manifold of dimension $n$ even. Then 
\begin{equation}\label{tractor4obstrn}
\fb_{n/2-2} W_{A_1A_2B_1B_2} =
K(n) X_{A_1}Z_{A_2}{}^aX_{B_1}Z_{B_2}{}^b \cB_{ab}.
\end{equation}
Here $K(n)$ is a known non-zero constant depending on $n$.  The tensor
$\cB_{ab} \in \ce_{(ab)_0}[2-n]$ is the Fefferman-Graham
obstruction tensor.  It is conformally invariant and natural.
\end{theorem}
\noindent We have $K(4)=-8$.  In dimensions $n\geq 6$, $K(n)$ is given
by $(n-4)k(n)$, where $k(n)$ is given in \nn{kn}, below.  Note that
$\fb_{n/2-2} W\in \ct^{2,2}[2-n]$. The theorem states that its
components vanish in all factors of the composition series
\nn{wcomp} for $\ct^{2,2}[2-n]$, except for the (injecting) factor
\mbox{$\ce_{1,1}[-2]=\ce^{1,1}[2-n]$}, and the term here is, up to scale, the obstruction.

{}From these theorems we have the following result.
\begin{theorem}\label{ceinv}
The obstruction tensor $\cB_{ab}$ vanishes on conformally Einstein
manifolds.
\end{theorem}
\noindent{\bf Proof:} A conformally Einstein manifold $M$ admits a
parallel standard tractor $\Pa$ (see \cite{GoNur}) such that $
\si:=\Pa^AX_A\neq 0$ is an Einstein scale.  It follows immediately
that $\Pa$ annihilates the tractor curvature $\Omega_{bc}{}^D{}_E$:
$\nd_c \Pa^D=0$ $\implies$
$\Omega_{bc}{}^D{}_E\Pa^E=[\nd_b,\nd_c]\Pa^D=0$. Also since $\Pa$ is
parallel then, viewing it as a multiplication operator, it is clear
that $[D,\Pa]=0$.  {}From \nn{WR} (see also \cite{CapGoamb}) we have $
W_{A_1A_2}{}^D{}_E = \frac{3}{n-2}D^{A_0}X_{A_0}
Z_{A_1}{}^bZ_{A_2}{}^c\Omega_{bc}{}^D{}_E$.  Thus $W_{BCDE}\Pa^E =0 $.

By Theorem \ref{powerslap} there is a formula for $\fb_{n/2-2}
W_{A_1A_2B_1B_2} $ which is polynomial in $\Box$, $D$, $W$,
$X$, $h$, and $h^{-1}$, and in this formula each of the indices
$A_1$, $A_2$, $B_1$, and $B_2$ appears on a $W$ tractor. On the other hand,
since $\Pa$ is parallel and of weight 0, it commutes with the operators
in this expression for $\fb_{n/2-2} W_{A_1A_2B_1B_2} $. Thus
\begin{equation}\label{killBoxW}
\Pa^{B_1} \fb_{n/2-2} W_{A_1A_2B_1B_2}=0
,
\end{equation}
since $ \Pa^{A}W_{ABCD}=0 $.

{}From \cite{GoNur} we have $\Pa^A=\frac{1}{n}D^A \si$.
Thus from the expression \nn{Dform} for
the tractor-D operator, we have the expression
$$
[\Pa^A]_g = \si Y^A -\frac{1}{n}\J \si X^A
$$ for $\Pa^A$ in terms of the (Einstein) metric
$g:=\si^{-2}\bg$. (Recall that if $\nabla$ is the Levi-Civita
connection determined by $g=\si^{-2}\bg$, then tautologically $\nabla
\si=0$.)  In particular, in this scale, we have
$
\Pa^AZ_A{}^a=0. 
$ 
Thus from Theorem \ref{main2} above,
$$
4(K(n))^{-1}Z^{A_2}{}_aZ^{B_2}{}_b\Pa^{A_1}\Pa^{B_1}\fb_{n/2-2} W_{A_1A_2B_1B_2}
=\si^2 \cB_{ab} .
$$
But from \nn{killBoxW}, the left-hand side vanishes, and hence
$\cB_{ab}=0$ on $M$. 
 \quad $\Box$
 
 Obtaining the obstruction tensor via a conformally invariant operator
 on a tractor field, as in Theorem \ref{main2}, enables us to relate
 it to other conformally invariant operators associated with the
 deformation complex, by ideas along the lines of the {\em curved
 translation principle} of Eastwood et alia \cite{Esrni}.  This is 
the idea behind
 Theorem \ref{obstrnthm}, which we are now ready to prove.  Related generalisations of the curved translation principle have been explored in depth in the
 setting of operators on differential forms \cite{BrGodeRham}.\\
\noindent{\bf Proof of Theorem \ref{obstrnthm}:} We first construct 
$\bB$ and prove \IT{iii}. Let $W^{2,2}$ denote the quotient of
$\bT^{2,2}[-2]$ by the subbundle which is the kernel of the bundle map
$\bT^{2,2}[-2]\to \bT^3\otimes \bT^3$ given by
$$
U_{A_2A_3B_2B_3}\mapsto X_{A_1}X_{B_1} U_{A_2A_3B_2B_3}
.
$$ We write $W_{2,2}$ for the subbundle of $\bT^{2,2}[2-n]$ consisting
of tractors which are annihilated by any contraction with $X$. We
write $\cW^{2,2}$ and $\cW_{2,2}$ for the section spaces of,
respectively, $W^{2,2}$ and $W_{2,2}$.  Note that complete
contractions between elements of $\bT^{2,2}[-2]$ and sections of
$\bT_{2,2}[2-n]$ take values in $E[-n]$. Hence there is a conformally
invariant pairing between $\ct^{2,2}[-2]$ and $\ct^{2,2}[2-n]$.
It is clear that the contractions between
elements of $\bT^{2,2}[-2]$ and sections of $\bT_{2,2}[2-n]$ induce a
well-defined bundle map $\langle \cdot ,\cdot \rangle:W^{2,2}\otimes
W_{2,2}\to E[-n]$, and so there is also a conformally invariant pairing
between $\cW^{2,2}$ and $\cW_{2,2}$.

Given a section $U_{ABCD}\in \ct^{2,2}[-2]$, let us write $[U_{ABCD}]$
for its image in the quotient space $\cW^{2,2}$.  {}From the tractor
composition series \nn{trcomp} (see also \nn{formtractorcomp} and the
discussion there), it follows easily that the space $\cW^{2,2}$ has a
composition series $\ce^{1,1}[2]\lpl \ce^{2,1}[2]\lpl \ce^{2,2}[2]$
and that the injection ${\rm I}:\ce^{2,2}[2] \to \cW^{2,2}$ is given
by
$$
u_{abcd}\mapsto [Z_A{}^aZ_B{}^bZ_C{}^cZ_D{}^d u_{abcd}].
$$

The differential operator $\bD:\cW^{2,2}\to \ct^{2,2}[-2]$ given by
$$
[U_{A_2A_3B_2B_3}] \mapsto \frac{9}{n(n-2)}\bY_{2,2} D^{A_1}D^{B_1}X_{A_1}X_{B_1}U_{A_2A_3B_2B_3}
$$
is clearly well-defined and conformally invariant. 
 Here $\bY_{2,2}$ is the bundle map which executes
the projection of $\bT^2[-1]\otimes \bT^2[-1]$ onto the direct summand
$\bT^{2,2}[-2]$.  We write $\bD^*$ for the formal adjoint of $\bD$.
This is a conformally invariant operator
$$
\bD^*: \ct^{2,2}[2-n]\to \cW_{2,2}.
$$
 On the other hand, from Theorem~\ref{powerslap}
there is a conformally invariant Laplacian type operator
$\fb_{n/2-2}:\ct^{2,2}[-2]\to \ct^{2,2}[2-n]$. Thus we have the
composition
$$
\bD^*\fb_{n/2-2}\bD : \cW^{2,2}\to \cW_{2,2}.
$$
The operator $\bB$ in the theorem is (up
to a constant multiple) simply the composition
$$
(\bD^*\fb_{n/2-2}\bD {\rm I} =:\bB): \ce^{2,2}[2]\to \cW_{2,2} .
$$
By construction this is natural and conformally invariant.

Now in a conformal scale, $(\bD {\rm I}(u))_{BCEF}$ is given explicitly by the
expression
\begin{equation}\label{DIform}
\begin{array}{c}
(n-4)\left((n-3)Z_B{}^bZ_C{}^cZ_E{}^eZ_F{}^f u_{bcef}\right. \\
\left. -2Z_B{}^bZ_C{}^cX_{[E}Z_{F]}{}^f \nd^e u_{efbc}   
-2 X_{[B}Z_{C]}{}^c Z_E{}^eZ_F{}^f \nd^b u_{bcef}  \right)\\
+ 4 X_{[B}Z_{C]}{}^c X_{[E} Z_{F]}{}^f(\nd^{(b}\nd^{e)}u_{bcef}+ (n-3)\Rho^{be} u_{bcef}). 
\end{array} 
\end{equation}
Thus from \nn{Wform} and a minor calculation,
\begin{equation}\label{recW}
\bD\big({\rm I}(C)\big)_{ABCD}= (n-3)W_{ABCD},
\end{equation}
 where $C$ is the Weyl curvature.

{}From Theorem~\ref{main2} and \nn{recW} we have
$$
(\fb_{n/2-2}\bD {\rm I} C)_{A_2A_3B_2B_3}=(n-3)K(n) 
X_{A_2}Z_{A_3}{}^aX_{B_2}Z_{B_3}{}^b\cB_{ab}
.
$$ That is, $\fb_{n/2-2}\bD {\rm I} C$ takes values in the factor
$\ce_{1,1}[-2]$ in the composition series for $\ct^{2,2}[2-n]$.  (Note
that this factor is a conformally invariant subspace.)  Now the formal
adjoint of the tractor-D operator is again the tractor-D operator
\cite{BrGoPacific}. So
$$
\begin{array}{lll}
\lefteqn{\bD^*X_{A_2}Z_{A_3}{}^aX_{B_2}Z_{B_3}{}^b\cB_{ab}=}&&
\\
&&
\frac{9}{n(n-2)}
X^{B_1}X^{A_1}D_{B_1}D_{A_1}
X_{A_2}Z_{A_3}{}^aX_{B_2}Z_{B_3}{}^b\cB_{ab}.
\end{array}
$$
But a short calculation using \nn{connids} and $\nn{Dform}$ shows
that this operation just returns $ 4(n-4)(n-3)
X_{A_2}Z_{A_3}{}^aX_{B_2}Z_{B_3}{}^b\cB_{ab}$, and this proves part
\IT{iii} of the theorem. All non-vanishing multiples can
be absorbed into the definition of $\bB$.

We treat now part \IT{i}. We need to show that ${\rm I}^*\bB$ has
order $n-4$ and is non-trivial. Since by construction there is a
universal natural expression for the operator $L$, it is sufficient to
establish this on the standard conformal sphere.  Recall that
$\fb_{n/2-2}$ has leading term $\Delta^{n/2-2}$.  Thus $\fb_{n/2-2}$
is elliptic (since the sphere has Riemannian signature).  {}From
\nn{DIform} it is clear that $\bD {\rm I}:\ce^{2,2}[2]\to
\ct^{2,2}[-2]$ is a differential splitting operator; there is a bundle
homomorphism ${\rm J}:\bT^{2,2}[-2]\to E^{2,2}[2]$ such that ${\rm
  J}\bD{\rm I}$ is the identity on $\ce^{2,2}[2]$. Thus on any
manifold, $\cR(\bD {\rm I}:\ce^{2,2}[2]\to \ct^{2,2}[-2])$ is
infinite-dimensional, and it follows immediately that $\fb_{n/2-2} \bD
{\rm I}$ is non-trivial on the standard conformal sphere. Now acting
on $\ce^{2,2}[2]$, $\fb_{n/2-2} \bD {\rm I}$ takes values in
$\ct^{2,2}[2-n]$.  The composition series for $\ct^{2,2}[2-n]$ is
given by \nn{wcomp} with $w=2-n$.  From this we see, for example, that
there is a canonical projection $\ct^{2,2}[2-n]\to
\ce^{1,1}[6-n]=\ce_{1,1}[2]$.
One can compose the operator $\fb_{n/2-2} \bD {\rm I}:\ce^{2,2}[2]\to
\ct^{2,2}[2-n]$ with this canonical projection.
 By construction, this is a conformally invariant
operator $\ce^{2,2}[2]\to \ce_{1,1}[2]$. On the other hand, from the
the classification of operators on conformally flat structures, as
discussed in Section \ref{def}, the only conformally invariant
operators on $\ce^{2,2}[2]$ taking values in irreducible bundles are as follows: there is an operator
$\ce^{2,2}[2]\to \ce^{3,2}[2]$ and an operator $\ce^{2,2}[2]\to
\ce_{2,2}[-2]$. {}From elementary weight considerations, we know the
latter has order $n-4$. Thus the composition described must be
trivial. \label{comps}  Continuing in this fashion and also using \nn{XYZtrans}, one
concludes that $\fb_{n/2-2} \bD {\rm I}$ takes values in the subspace
$\cW_{2,2}=\ce_{2,2}[-2]\lpl \ce_{2,1}[-2]\lpl \ce_{1,1}[-2]$, and the
composition of $\fb_{n/2-2} \bD {\rm I}$ with projection to
$\ce_{2,2}[-2]$ is necessarily non-trivial. This composition is thus
the unique (up to scale) conformally invariant operator between these
bundles (on the conformal sphere). We are now done as follows. On the
one hand, ${\rm I}^*\bD^*$ is the formal adjoint of a splitting
operator for $\ce^{2,2}[2]$ and therefore acts as a multiple of the
identity on the component $\ce_{2,2}[-2]$. On the other hand, ${\rm
I}^*\bD^*$ must annihilate the components $\ce_{2,1}[-2]$ and
$\ce_{1,1}[-2]$, since these have higher total order than the target
bundle for the composition (viz.\ $\ce_{2,2}[-2]$) and a natural
differential operator cannot lower order.

Finally, we consider \IT{ii}. Let us first consider the case of a flat
Riemannian or pseudo-Riemannian structure.  So all curvature will
vanish, until we note otherwise.  Let $F_{a_1a_2}$ denote a
2-form. Then
$\Delta^{n/2-2}D^{A_0}X_{A_0}Z_{A_1}{}^{a_1}Z_{A_2}{}^{a_2}F_{a_1a_2}
$ is well understood as a special case of the results in
\cite{BrGodeRham}.  (See Proposition 4.6.)  The non-zero components of
this have values in a subbundle of $\ct^2[2-n]$ with composition
series $\ce_2\oplus\ce_1$. These components are (up to an overall
non-zero constant multiple)
$$
\left(\begin{array}{l}(\d d)^{n/2-2}F\\ 
a \d(d\d)^{n/2-2} F \end{array} \right) ,
$$
where $d$ is the exterior derivative, $ \d $ its formal adjoint,
and $a$ is a non-zero constant. Composing these components on the left
with $(\d d~,~\frac{1}{a}d)$ yields $(\d d +d\d)^{n/2-1}F=(-1)^{n/2-1}
\Delta^{n/2-1} F$. Now on flat structures we have the identity 
$$
(n-2)\big(\bD {\rm I} (u)\big)_{A_1A_2B_1B_2}= 3
D^{A_0}X_{A_0}Z_{A_1}{}^{a_1}Z_{A_2}{}^{a_2}U_{a_1a_2B_1B_2},
$$
where $u\in \ce^{2,2}[2]$ and $ U_{a_1a_2B_1B_2}$ is the
conformally invariant form-tractor given in scale by letting $w$
equal 2 in the formula
\begin{equation}\label{trext}
U_{a_1a_2B_1B_2}= (n+w-5)Z_{B_1}{}^{b_1}Z_{B_2}{}^{b_2} u_{a_1a_2b_1b_2}+
2X_{B_1}Z_{B_2}{}^{b_2} \nd^{b_1} u_{a_1a_2b_2b_1}.
\end{equation}
Thus by viewing $U$ as a 2-form with values in a tractor bundle and
replacing $\nd$, $d$, and $\d$ with their tractor connection coupled
variants in the argument above, we conclude that there is an operator
$\overline{\bA}$ such that
$$
\begin{aligned}
\overline{\bA}\Delta^{n/2-2}\bD {\rm I} (u) &=
\frac{3}{n-2}
\overline{\bA}\Delta^{n/2-2}
D^{A_0}X_{A_0}Z_{A_1}{}^{a_1}Z_{A_2}{}^{a_2}U_{a_1a_2B_1B_2}\\
&
=\Delta^{n/2-1} U.
\end{aligned}
$$ We continue with similar considerations, except that now we view
$u^{a_1}{}_{a_2b_1b_2}$ as a 2-form on the $b_1b_2$ index pair that
takes values in End$(TM)$. If $F$ now indicates a 2-form of weight
$w'$, then we have
$$
\begin{aligned}
\bK(F) :&=
\frac{3}{n+2w'-2}D^{A_0}X_{A_0}Z_{A_1}{}^{a_1}Z_{A_2}{}^{a_2}F_{a_1a_2}
\\
&= (n+w'-4)Z_{A_1}{}^{a_1}Z_{A_2}{}^{a_2}F_{a_1a_2}+
2X_{A_1}Z_{A_2}{}^{a_2}\nd^{a_1}F_{a_2a_1} .
\end{aligned}
$$ So if, in particular, $w' =1$, then the formula on the right-hand
side agrees with \nn{trext} (with $w=2$). In formally calculating
$\Delta^{n/2-1} U^{a_1}{}_{a_2B_1B_2}$ using the identities
\nn{connids} and the Leibniz rule to obtain a formula polynomial in
$u$, $\nd$, the metric $g$, its inverse, and the projectors $X$, $Y$, and
$Z$, we may ignore the $a_1$ and $a_2$. Their contribution is buried
in the meaning of the Levi-Civita connection $\nd$.  Now for a 2-form
$F$ of weight $1$, we have that on flat structures, $ \Delta^{n/2-1}
\bK(F) $ takes values in $\ce_2[-1]\lpl \ce_1[-1]$ and has the
form \cite{BrGodeRham}
$$
\left(\begin{array}{c}\big((3-n) (\d d)^{n/2-1} + (d\d)^{n/2-1} \big) F\\ 
{\ast} \end{array} \right)
$$
there, up to an overall non-zero multiple. Here $\ast$ indicates some
term, the details of which will not concern us.  We note that the top
expression gives an elliptic operator on $F$; we may act on this with
the operator $\d d+(3-n)d\d$ to yield
$(3-n)(-1)^{n/2}\Delta^{n/2}F$.  Thus there is a linear differential
operator $\overline{\bA}_2$ such that
$$
\overline{\bA}_2 \Delta^{n/2-1}U = \Delta^{n/2}u.
$$

Combining these observations, we see that there is a linear differential
 operator $\overline{\bA}_3$ such that 
$$
\overline{\bA}_3 \Delta^{n/2-2}\bD {\rm I} (u) =\Delta^{n/2} u .
$$
Finally, we note that one can easily verify directly that $\bD^*$ is
differentially invertible as a graded differential operator on the
subspace $\cW_{2,2}$. (That is, its inverse is also a graded
differential operator. The point is that in terms of a splitting of
$\cW_{2,2}$ determined by a choice of conformal scale, a
straightforward calculation shows that $\bD^*$ takes the form
$(u,v,w)\mapsto (ku,\ell v+ \delta\cdot u,m w +\delta\cdot v
+\delta\cdot\delta\cdot u )$, where $k$,
$\ell$, and $m$ are non-zero integers, $\delta\cdot$ indicates a
divergence operator, and $\delta\cdot \delta\cdot$ a double divergence
operator.)  Thus with $\overline{\bB}$ defined to be the necessary
multiple of $ \overline{\bA}_3 (\bD^*)^{-1}$, we have the result
$\IT{ii}$ for flat structures. But now the result follows in general,
since moving to curved structures yields the same formal calculation,
except that at each stage the differential operators concerned may
have additional lower order terms involving curvature. It is easily
checked that these terms can only yield terms of order lower than $n$
in the final calculation of $\overline{\bB}\bB$.  \quad $\Box$\\

\noindent{\bf Proof of Proposition \ref{factoring}:} We treat $\sf L$ first.
We already have $\sf L= {\sf Bi}^*{\sf Bi}$ in  dimension 6, and so  we shall 
assume that $n\geq 8$. Let us denote by 
$$
\bU:\ce^{2,2}[2]\to \ce^{2}\otimes\ct^2
$$
the conformally invariant operator given by \nn{trext}. We write
$d^{\nd}$ for the tractor connection coupled exterior derivative and
$\d^\nd$ for its formal adjoint. Thus for example for $U\in
\ce^{2}\otimes\ct^2 $ we have $(d^{\nd}U)_{a_0a_1a_2B_1B_2}=
3\nd_{a_0}U_{a_1a_2B_1B_2} $.  It is straightforward using
\nn{connids} to verify that on conformally flat structures, the composition 
 $d^\nd \bU $ can be re-expressed in the form $\bM
{\sf Bi} $, where $\bM :\ce^{3,2}[2]\to \ce^{3}\otimes \ct$ is a
conformally invariant first-order differential splitting operator.

There are conformally invariant formally self-adjoint operators
$L_{k}:\ce^k \to \ce_k$, $0\leq k\leq n$, with leading term $(\d
d)^k$.  These are the ``long operators'' for the de Rham complex given 
 in \cite{BrGodeRham}. It
is shown there that there are natural linear differential operators
$Q_{k+1}$ such that $L_k=\d Q_{k+1} d$.

Now suppose we are on a contractible (conformally flat) manifold. This
suffices for our present purposes. Then the tractor bundle is flat and
trivial.  It follows that there are conformally invariant and formally
self-adjoint tractor-coupled variants of the $L_k$,
$$
L_k^{\nd}: \ce^k\otimes\ct^2 \to \ce_k\otimes \ct^2.
$$
These are obtained by formally replacing, in the natural formulae
for  $L_k=\d Q_{k+1} d$, each instance of $d$, $\d$, and the
Levi-Civita connection with, respectively, $d^\nd$, $\d^\nd$, and the
Levi-Civita tractor-coupled connection. By construction the result has a
factorisation $L_k^{\nd} =\d^\nd Q_{k+1}^\nd d^\nd$ for some
differential operator $Q_{k+1}^\nd$. 

Observe that by composition, we have a formally self-adjoint
conformally invariant operator $\bU^* L_2^{\nd} \bU:\ce^{2,2}[2]\to
\ce_{2,2}[-2]$, where $\bU^*$ is the formal adjoint of $\bU$. We will
re-express this. By taking formal adjoints, we have $\bU^* \d^\nd= {\sf
Bi}^* \bM^*$ from $d^\nd \bU =\bM {\sf Bi} $.  Thus we obtain an
operator
$$
{\sf Bi}^*  \bM^* Q_{3}^\nd \bM {\sf Bi} : \ce^{2,2}[2]\to \ce_{2,2}[-2] .
$$ The result follows from the uniqueness of $\sf L$, provided the
displayed operator is non-trivial. It is clearly sufficient to
establish this for Riemannian signature structures and at a flat
metric within the conformal class. We use the alternative expression
$\bU^* \d^\nd Q_{3}^\nd d^\nd \bU = {\sf Bi}^* \bM^* Q_{3}^\nd \bM
{\sf Bi}$.  On flat structures, $Q_3=(d\d)^{n/2-3}$, and so
$Q_3^\nd=(d^\nd \d^\nd)^{n/2-3}$. Since the tractor connection is
flat, it follows that for $u\in \ce^{2,2}[2]$ of compact support, $
\d^{\nd} (d^\nd \d^\nd)^{n/2-3}d^\nd \bU(u) $ vanishes if and only if
$d^\nd \bU(u) $ vanishes. (Suppose $d^\nd \bU(u)\neq 0$. Then there
exists a parallel $T\in \ct^2$ such that
$T^{B_1B_2}(d^{\nd}\bU(u))_{a_0a_1a_2B_1B_2}\neq 0 $, i.e., $d f\neq
0$, where $f_{a_1a_2}:= T^{B_1B_2}\bU(u)_{a_1a_2B_1B_2}$. But, on the
other hand, if $0= T^{B_1B_2}(\d^{\nd} (d^\nd \d^\nd)^{n/2-3}d^\nd
\bU(u))_{a_0a_1a_2B_1B_2}$, then $(\d d)^{n/2-2}f=0$ $\Rightarrow$
$df=0$.)  This is equivalent to $\bM {\sf Bi}(u) $ vanishing.  Since
$\bM$ is a differential splitting operator, this in turn is equivalent
to ${\sf Bi}(u)=0$. Thus the composition $ \d^\nd Q_{3}^\nd d^\nd \bU
:\ce^{2,2}[2]\to \ce_2\otimes \ct^2 $ is non-trivial.  Now it is
easily verified that $\ce_{2,2}[-2]$ turns up with multiplicity 1 in
the composition series for $ \ce_2\otimes \ct^2 $. It follows, by an
exact analogue of the argument used on page \pageref{comps}
that $ \d^\nd
Q_{3}^\nd d^\nd \bU $ takes values in, and only in, $\ce_{2,2}[-2]$
and composition factors of higher total order. Thus on the range of
this operator, $\bU^*$ acts as a non-zero multiple of the projection
to the component $\ce_{2,2}[-2]$. (Recall that $\bU^*$ is the formal
adjoint of a differential splitting operator $\bU:\ce^{2,2}[2]\to
\ce^{2}\otimes\ct^2$, and so it must act as a non-zero multiple of the
identity on the component $\ce_{2,2}[-2]$. On the other hand, it is
differential, so it cannot lower total order.)

Now we consider the situation for $\sf B$. 
We require a conformally invariant differential
splitting operator $\Gamma: \ce^{1,1}[2]\to \ce^1\otimes \ct^2$ that
will play a role, in this case, analogous to the role of $\bU$ above.
This is easily constructed explicitly and directly (and can be
obtained from a composition of the related operators in Section~5.1 of
\cite{BrGoopava}), and so we omit the details. Since $\Gamma$ has
values in a weight zero adjoint tractor-valued bundle of 1-forms it is
clear that the composition $d^\nd \Gamma $ is conformally
invariant. This is easily verified non-trivial.
 On the other hand, in terms of a metric $g$, the tractor
curvature is given by
$$
Z_{B_1}{}^{b_1}Z_{B_2}{}^{b_2} C_{a_1a_2b_1b_2}+
\frac{2}{n-3}X_{B_1}Z_{B_2}{}^{b_2} \nd^{b_1} C_{a_1a_2b_2b_1}.
$$
Thus the linearisation, at a conformally flat metric $g_0$, of the
tractor curvature is $\frac{1}{n-3}\bU {\sf C} $. This is manifestly
non-trivial, and so, via arguments as used several times already
concerning the uniqueness of irreducible conformally invariant
operators, it is straightforward to verify that this operator must
agree with $d^\nd \Gamma $ (on conformally flat structures), at least up to scale.  We set the scale
of $\Gamma$ so that $ d^\nd \Gamma = \frac{1}{n-3}\bU {\sf C} $.  On flat
manifolds, $Q_2=(d \d)^{n/2-2}$, and so by almost the same argument as
for $\sf L$, we conclude that on conformally flat manifolds, the
formally self-adjoint conformally invariant operator ${\sf C}^*\bU^*
Q^\nd_2 \bU {\sf C}$ is non-trivial.  \quad $\Box$

The next theorem shows that for $n$ even, if the ambient curvature is
formally Ricci flat to $O(Q^{n/2-1})$, then a tensor part of the
coefficient of $Q^{n/2-1}$ is a natural conformal invariant of the
underlying manifold and so is an obstruction to finding an ambient
metric which is Ricci flat to higher order. For our purposes, the main
point is that this is achieved by \nn{RtoB}, which recovers this
obstruction via a tangential operator acting on the ambient curvature.
\begin{theorem} \label{main1}
 For a conformal manifold $M$ of even dimension $n$, let $\h$ be an
associated ambient metric satisfying
$\Ric(\h)=Q^{n/2-1}\aB$.  Then we have\\
\IT{i} $\aB|_\cq$ is equivalent to a tractor
  $B_{AB}\in \ce_{(AB)_0}[-n]$ such that $X^A B_{AB}=0$.\\
\IT{ii} The weighted tensor $Z^A{}_aZ^B{}_bB_{AB}=:\cB_{ab}$
is a section of $ \ce_{(ab)_0}[2-n]$.
\IT{iii}  For $n\geq 6$, we have 
\begin{equation}\label{RtoB}
\afl^{n/2-2} \aR_{A_1A_2B_1B_2} 
= k(n) \X_{A_1}\X_{B_1}\aB_{A_2 B_2} +O(Q) ,
\end{equation}
where $k(n)$ is the dimension dependent non-zero constant given
above.\\ 
In dimension 4, 
$$
3\D^{A_0}\X_{A_0}\aR_{A_1A_2B_1B_2}=
 16 \X_{A_2}\X_{B_1}\aB_{A_1B_2} +O(Q)
.
$$
\IT{iv} The tensor $\cB_{ab} $ is divergence-free.
\\
\IT{v} The weighted tensor $\cB_{ab}$ is a non-trivial natural
conformal invariant of the form
$\Delta^{n/2-2}\nd^c\nd^d C_{cadb} + \LOT =(n-3)\Delta^{n/2-2} (\Delta
\Rho_{ab}-\nd_a\nd_b\J) +\LOT $
(up to a constant multiple),
and so is an obstruction to finding an
ambient metric which is Ricci flat modulo $O(Q^{n/2})$.
\end{theorem}
\noindent{\bf Remarks:} 1. The statement of the
theorem up to the definition of $\cB_{ab}$ in \IT{ii} is a
characterisation of the Fefferman-Graham obstruction tensor
\cite{Grpriv}.  See also \cite{FeffGrinprogess}. (This gives a complete
obstruction to the ambient metric in the sense that if this vanishes,
then the ambient construction may be continued to all orders
\cite{FGast}.)
Hence $\cB_{ab}$ is the usual obstruction tensor, as claimed in
Theorem~\ref{main2}.
Thus part \IT{iii}, above, gives a new
ambient formula for the Fefferman-Graham obstruction tensor.\\
2. {}From \nn{abracketQ} it follows easily that $\cB_{ab}$ may be equally
 viewed as an obstruction to obtaining an ambient metric which is
 harmonic for $\afl$ in the sense that $\afl \aR$ vanishes to all
 orders. See also the remark on page~\pageref{harmony}.\\
3. It should be pointed out that  
\begin{equation}\label{altRtoB}
\al^{n/2-3}\afl \aR_{A_1A_2B_1B_2} 
= k(n) \X_{A_1}\X_{B_1}\aB_{A_2 B_2} +O(Q)
\end{equation}
is an alternative ambient formula for the obstruction, and we could
replace the $\al^{n/2-3}$ by $\al^{n/2-3}_\alpha$ in this
formula. \quad \endrk

\noindent{\bf Proof of Theorems \ref{main2} and \ref{main1}:} As
above, we write $\aRic$ for $\Ric(\h)$.
It is immediate that $\aB$ is symmetric and homogeneous of weight
$-n$. Also from \nn{XRQ} it follows that $\X^A\aB_{AB}=0$.  So
$\aB|_\cq$ is equivalent to a tractor field $B_{AB}\in \ct_{(AB)}[-n]$
satisfying $X^A B_{AB}=0$. {}From this last equality and \nn{XYZtrans},
it is clear that $\cB_{ab}$ is conformally invariant, while from the
weight and symmetry of $B_{AB}$, it follows that $\cB_{ab}\in
\ce_{(ab)}[2-n]$. For parts \IT{i} and \IT{ii}, it remains to show
that both $B_{AB}$ and $ \cB_{ab}$ are trace-free.

First we consider the case $n\neq 4$.  Note that since $\nda_A
Q=2\X_A$, we have
\begin{equation}\label{skewRic}
  \nda_{A_1}\nda_{B_1}\aRic_{A_2 B_2} =
  (n-2)(n-4)Q^{n/2-3}\X_{A_1}\X_{B_1}\aB_{A_2 B_2}+O(Q^{n/2-2})
  .
\end{equation}
{}From \nn{semihar} and \nn{abracketQ} together with a short
computation, it follows that
$$%
\afl^{n/2-2}\aR_{A_1A_2B_1B_2}  
= k(n) \X_{A_1}\X_{B_1}\aB_{A_2 B_2} +O(Q) ,
$$
as claimed in \IT{iii}, where  
\begin{equation}\label{kn}
  k(n)=
  (n-2)(n-4)(-1)^{n/2-3}2^{n-4}\left((n/2-3)!\rule{0mm}{4mm}\right)^2
  .
\end{equation}
(Note that \nn{semihar} and \nn{abracketQ} also give the alternative
formula in Remark 3, above.)

Since $(n-4)\aR|_\cq$ is equivalent to the tractor field $W$, it
follows from Proposition~\ref{powerslaptan}, below, that
$(n-4)\afl^{n/2-2} \aR|_\cq$ descends to the natural tractor field
$\fb_{n/2-2}W$.  On the other hand, using
$\delta_{B}{}^A=X^AY_B+Y^AX_B+ Z^A{}_{a}Z_{B}{}^b\delta_b{}^a$ and the
fact that $X^AB_{AB}=0$, we see that
\begin{equation}\label{Lab17May4a}
X_{A_1}X_{B_1}B_{A_2B_2}=
X_{A_1}Z_{A_2}{}^aX_{B_1}Z_{B_2}{}^b\cB_{ab}.
\end{equation}%
It therefore follows that $ \X_{A_1}\X_{B_1}\aB_{A_2 B_2}|_\cq$ is
equivalent to the tractor field $X_{A_1}Z_{A_2}{}^a X_{B_1}Z_{B_2}{}^b
\cB_{ab}$.  This establishes \nn{tractor4obstrn} of
Theorem~\ref{main2}.

Since left-hand side of \nn{tractor4obstrn} is natural, it follows
that $X_{A_1}X_{B_1}B_{A_2B_2}$ is natural.  Hence
$$
\cB_{ab} =Z^{A_2}{}_aZ^{B_2}{}_bB_{A_2 B_2} 
=4Y^{A_1}Y^{B_1}Z^{A_2}{}_aZ^{B_2}{}_b X_{A_1}X_{B_1}B_{A_2 B_2}
$$
is likewise natural, as claimed in \IT{v} and Theorem~\ref{main2}.

Next we show that $B_{AB}$ and $\cB_{ab}$ are trace-free.  According
to Theorem~\ref{powerslap}, the operators $\fb_m$ preserve tensor
type.  Since $W_{A_1A_2B_1B_2}$ is trace-free, it follows that
$\fb_{n/2-2}W_{A_1A_2B_1B_2}$ is completely trace-free.  Thus from
\nn{tractor4obstrn} and \nn{Lab17May4a} it follows that $
h^{A_1B_2}X_{A_1}X_{B_1}B_{B_2 A_2} =0.  $ Since $B_{CD}$ is symmetric
and $X^A B_{AB}=0$, it follows that $h^{AB}B_{AB}=0$ as claimed.  Now
using \nn{innerprodt} and once again the fact that
$X^AB_{AB}=0$, we see that that $\bg^{ab}\cB_{ab}=0$.

We must obtain the corresponding results in dimension 4. First observe
that in any dimension,
\[
\begin{array}{lll}
\lefteqn{3\D^{A_0}\X_{A_0}\aR_{A_1A_2B_1B_2}=}&&\vspace{3mm}
\\
&&
(n-2)[(n-4)\aR_{A_1A_2 B_1B_2}+2\X_{A_1}\nda^{C}\aR_{A_2CB_1B_2}] +O(Q),
\end{array}
\]
by \nn{alternativeD} and \nn{DAboldXA}.  {}From the contracted Bianchi
identity, we have, for $n=4$,
$$
\begin{array}{lll}
\lefteqn{3\D^{A_0}\X_{A_0}\aR_{A_1A_2B_1B_2}=}&&\vspace{3mm}
\\
&&
8 \X_{A_2}\nda_{B_1}\aRic_{A_1B_2} +O(Q)
= 16 \X_{A_2}\X_{B_1}\aB_{A_1B_2} +O(Q)
.
\end{array}
$$ Relating $W$ to the left-hand side via \nn{WR}, we conclude that in
dimension 4,
$$
W_{A_1A_2B_1B_2}= -8 X_{A_1}X_{B_1} B_{A_2 B_2},
$$
and comparing this with the  formula \nn{Wform} for $W$ above, we have
\[
-2 X_{A_1}X_{B_1} B_{A_2 B_2}= X_{A_1}Z_{A_2}{}^a X_{B_1} Z_{B_2}{}^b
B_{ab}.
\]
Thus $\cB_{ab}$ is a scalar multiple of the Bach tensor,
$
\cB_{ab} = -\frac{1}{2}B_{ab} ,
$ which is natural and trace-free, by \nn{one}. We note also that
since $W$ is trace-free and $X^AB_{AB}=0$, it follows that $B_{AB}$ is
trace-free.

It is well known (and easily verified) that the Bach tensor
in dimension 4 is divergence-free. For \IT{iv} we need the analogous
result in other dimensions. First note that a short calculation, which
uses the formula \nn{Dform} for the tractor-D operator and the
identities \nn{connids} for the connection, shows that
$$
2D^{A_1}X_{A_1}Z_{A_2}{}^a X_{B_1}Z_{B_2}{}^b
\cB_{ab}=(n-4)X_{A_2} X_{B_1}Z_{B_2}{}^b \nd^a \cB_{ab}. 
$$
So, in dimensions other than 4, it follows that $D^{A_1}X_{A_1} X_{B_1}
B_{A_2B_2}$, and equivalently $(\D^{A_1}\X_{A_1} \X_{B_1} \aB_{A_2B_2
})|_\cq$, vanish if and only if $ \nd^a \cB_{ab}=0$. We
calculate $\D^{A_1}\X_{A_1} \X_{B_1} \aB_{A_2B_2 }$ on the ambient
manifold.  By \nn{altRtoB}, this is
$$
 \D^{A_1}\al^{n/2-3} \afl \aR_{A_1A_2B_1B_2}+O(Q),
$$
up to a non-zero multiple, since $\D$ acts tangentially.  We ignore
terms $O(Q)$ for much of the remainder of this calculation.  The above
display expands to
$$ (4-n)\nda^{A_1}\al^{n/2-3}\afl \aR_{A_1A_2B_1B_2}-\al\X^{A_1}
\al^{n/2-3}\afl \aR_{A_1A_2B_1B_2}.
$$
{}From \nn{XRQ} and \nn{abracketalX} we obtain
\begin{equation}\label{DXXB}
\begin{array}{l}
(4-n)[\nda^{A_1},\al]\al^{n/2-4}\afl \aR_{A_1A_2B_1B_2 }
\vspace{1mm}
\\
+(6-n) \al [\nda^{A_1},\al]\al^{n/2-5}\afl \aR_{A_1A_2B_1B_2}
\vspace{1mm}\\
+\cdots - 4 \al^{n/2-4} [\nda^{A_1},\al]\afl \aR_{A_1A_2B_1B_2}
\vspace{1mm}
\\
-2 \al^{n/2-3}\big([\nda^{A_1},\al] \aR_{A_1A_2B_1B_2}
+\frac{1}{2}\nda^{A_1}(\aR\sharp\sharp\aR_{A_1A_2B_1B_2})\big)
,
\end{array}
\end{equation}
after some re-organisation. It remains only to observe that all
the terms in this sum are $O(Q)$. First we note that from \nn{semihar} and  
\nn{skewRic}, it is clear that 
$$
\afl \aR_{A_1A_2B_1B_2}= K Q^{n/2-3}\X_{A_1}\X_{B_1}\aB_{A_2 B_2}+O(Q^{n/2-2}) ,
$$ for some constant $K$.  Thus by \nn{abracketQ}, each term
$$
\al^k [\nda^{A_1},\al]\al^{\ell}\afl \aR_{A_1A_2B_1B_2},
\quad k+\ell=n/2-4,
$$
is some number times
\begin{equation}\label{Lab06Aug4a}
\al^k[\nda^{A_1},\al]Q^{n/2-3-\ell}\X_{A_1}\X_{B_1}\aB_{A_2
B_2}+O(Q)
,
\end{equation}
since $[\nda^{A_1},\al]$ is a first-order operator.  Now consider the
identity obtained from \nn{L17May01c} by including the
$\BigO(Q^{n/2-2})$ terms omitted from the display in \nn{L17May01c}.
{}From this identity, from \nn{XRQ}, and from the fact that $\nda
Q=2\X$, it follows that
$[[\nda^{A_1},\al],Q]=0$ identically on the ambient manifold.  Thus
\nn{Lab06Aug4a} is $\BigO(Q)$.

Now consider the last term in \nn{DXXB}.  By direct calculation, we
have
$$
[\nda^{A_1},\al] \aR_{A_1A_2B_1B_2}=
-\frac{1}{2}\nda^{A_1}(\aR\sharp\sharp\aR_{A_1A_2B_1B_2})+O(Q^{n/2-2}),
$$
and so 
$$
\al^{n/2-3}\big([\nda^{A_1},\al] \aR_{A_1A_2B_1B_2} 
+\frac{1}{2}\nda^{A_1}(\aR\sharp\sharp\aR_{A_1A_2B_1B_2})\big)=O(Q)
$$
as required. 

Finally, we must show that in general $\cB_{ab}$ is non-trivial. 
Up to scale, $\cB_{ab}$ is given by
$$
4Y^{A_1}Y^{B_1}Z^{A_2}_aZ^{B_2}{}_b \fb_{n/2-2}W_{A_1A_2B_1B_2} .
$$
{}From \nn{Wform} and \nn{bi1}, it is clear that
$4Y^{A_1}Y^{B_1}Z^{A_2}_aZ^{B_2}{}_b W_{A_1A_2B_1B_2} $ is, at leading
order, a non-zero multiple of $\nd^d\nd^c C_{cadb}$.  Using the fact
that $\fb_{n/2-2} $ has leading term $\Delta^{n/2-2}$, and then
\nn{connids} to verify that the commutator of $\Delta^{n/2-2}$ with
$4Y^{A_1}Y^{B_1}Z^{A_2}_aZ^{B_2}{}_b $ generates only lower order
terms, we conclude that
$$ 
\cB_{ab} = \ell(n) \Delta^{n/2-2}\nd^d\nd^c C_{cadb} + \LOT
,
$$
where $\ell(n)$ is a non-zero constant.  Given the form of the leading
term, an elementary exercise shows that this natural tensor cannot
vanish in general.  \quad $\Box$

\subsection{Conformal Laplacian operators on tractor fields} \label{lapsect}

It remains to prove Theorem~\ref{powerslap}.  Our strategy is to first
define the operators $\Box_m^\alpha$, which we do via powers of the
ambient Laplacian $\ala$ in Proposition~\ref{powerslaptan}, and then
rewrite each such power as a combination of compositions of low order
tangential operators, each of which has an immediate interpretation as
an operator on a tractor bundle.  This leads to a simple algorithm for
rewriting any operator of this form in terms of basic tractor
operators using only the existence of an ambient metric.  Two of the
key tools are Theorem~\ref{curv}, which explains how ambient
derivatives of the ambient curvature can be re-expressed in terms of
low order tangential operators, and Proposition~\ref{Lab13July3g},
which describes harmonic extensions of tensor fields along $\cq$.

Almost all of the subsequent discussion concerns the ambient manifold
$\tilde{M}$ with metric as discussed in Section \ref{ambient}.
Occasionally we pause to interpret results on the underlying conformal
manifold $M$.

In the generic $n$-even case, some identities, such as (\ref{aflaR})
and (\ref{L17May01c}), hold to only finite order in $Q$.  In many
proofs, we will apply the operators $\nda$ and $\al$ to both sides of
an identity, and this will reduce the order to which the identity
holds.  Thus we must keep track of the number of times that we apply
$\nda$ and $\al$.  In odd dimensions and in the conformally flat case,
this is unnecessary, since the identities hold to all orders.  For
simplicity, many of the proofs that follow explicitly treat only the
generic $n$-even case, since the proofs in the other cases are
essentially the same, except for the fact that they do not require the
operator counts.  In addition, we have stated some of the results
themselves in the generic $n$-even case only.  All results hold as
stated.  Propositions~\ref{lab23June3a} and \ref{lab23June3c},
Theorem~\ref{curv}, and Lemma~\ref{arranging} also hold in general;
they hold to all orders in both the odd-dimensional case and the
conformally flat case, and in these cases the upper bounds stated in
the \textit{hypotheses} of the results no longer apply.

We will often use abbreviated notations.
We may abbreviate \nn{L17May01c} by writing $[\al,\nda]V=\sum\aR\nda
V$.  It is easily verified that \nn{L17May01c} generalises to
\begin{equation}\label{abracketnda}
[\ala,\nda]V =
{\textstyle\sum}\aR\nda V
+ \alpha{\textstyle\sum}(\nda \aR)V,
\end{equation}
which also holds modulo $O(Q^{n/2-2})$ in even dimensions and to
infinite order in odd dimensions.  
For example, let $V$ be any symmetric ambient 2-tensor.  In this case
\nn{abracketnda} stands for
$$
\begin{array}{rll}
\lefteqn{[\ala,\nda_A]V_{BC}=}&&
\\
&&
  2(\alpha-1)\aR_A{}^P{}_B{}^Q\nda_P V_{QC}
   +2(\alpha-1)\aR_A{}^P{}_C{}^Q\nda_P V_{BQ}
\\
&&
   -2\alpha(\nda_A\aR_{B}{}^P{}_C{}^Q)V_{PQ}
   ,
\end{array}
$$
which holds to the appropriate order.  If the $V$ on the left-hand
side of \nn{abracketnda} has any free indices, then in every term of
the right-hand side of \nn{abracketnda}, each such index either
remains attached to $V$ in its original position or moves onto an
$\aR$.  Some of the proofs in Section~\ref{FormulaeSect} will use this
fact, which follows immediately from \nn{L17May01c} and the definition
of $\aR\hash\hash$.  The expressions we treat will often involve
iterations of operators.  To indicate how many operators we are
composing in such an iteration, we will use exponents.  For example,
we might indicate $\nda_A\nda_B\aR_{CDEF}$ by writing $\nda^2\aR$.
We will often use the symbol $\cp$ to denote a
partial contraction polynomial.  The same symbol $\cp$ may denote
different polynomials in different parts of a given discussion.

We often use the identities \nn{delQ} and $\nda Q=2\X$ without
explicit mention.

The proof of Theorem~\ref{powerslap} begins with the development of a
useful ambient calculus.  This involves a sequence of results.
\begin{proposition}
\label{lab23June3a}
Suppose that $n$ is even and $M$ is generic.  Let an integer $\ell$ be
given, and suppose that $0\leq\ell\leq\frac{n}{2}-4$.  Then on the
ambient manifold,
\begin{equation}\label{lab23June3b}
\al\nda^\ell\aR=
\sum(\nda^p\aR)(\nda^q\aR)+\BigO(Q^{n/2-3-\ell}),
\end{equation}
where $p+q=\ell$.  If the $\aR$ on the left-hand side of
(\ref{lab23June3b}) has any free indices, then for every term in the
summation, these indices appear on an $\aR$ (as opposed to a $\nda$).
\end{proposition}
\noindent{\bf Proof:} We use induction.  The case $\ell=0$ follows
from \nn{aflaR}.  Suppose next that $0\leq m\leq \frac{n}{2}-5$ and
that the result holds for $\ell=m$.  {}From this assumption and
(\ref{L17May01c}), we have
\[
\begin{array}{l}
  \al\nda^{m+1}\aR=
  \nda\al\nda^m\aR
  +\sum\aR(\nda^{m+1}\aR)
  +\BigO(Q^{n/2-2})=
\vspace{2mm}
\\
\nda\left(\sum(\nda^p\aR)(\nda^q\aR)+
  \BigO(Q^{n/2-3-m})\rule{0mm}{4mm}\right)
  +\sum\aR(\nda^{m+1}\aR)
  +\BigO(Q^{n/2-2})
  =
\\
\sum(\nda^s\aR)(\nda^t\aR)+\BigO(
\rule{0mm}{4mm}
Q^{n/2-3-(m+1)}
).
\rule{0mm}{6mm}
\end{array}
\]
Here $p+q=m$ and $s+t=m+1$.  The use of the inductive
assumption and \nn{L17May01c} never moves a free index from an $\aR$
onto a $\nda$.  $\Box$

\begin{proposition}
\label{lab23June3c}
Suppose that $n$ is even and $M$ is generic.  Let an integer $\ell$ be given,
and suppose that $0\leq\ell\leq\frac{n}{2}-3$.  Then
\begin{equation}\label{lab23June3d}
\al^\ell\aR=\sum(\nda^{v_1}\aR)\cdots(\nda^{v_j}\aR)+
\BigO(Q^{n/2-2-\ell}).
\end{equation}
In (\ref{lab23June3d}), the number of factors in a term may vary from
term to term, but in any case, $v_i\leq\ell$ for $1\leq i\leq j$.  If
$A$, $B$, $C$, and $D$ denote the indices of the $\aR$ on the
left-hand side of (\ref{lab23June3d}), then for each term in the sum,
these indices are on an $\aR$.
\end{proposition}
\noindent{\bf Proof:} We again use induction.  Suppose that $0\leq
m\leq\frac{n}{2}-4$ and that the result holds for $\ell=m$.  Then
\begin{equation}\label{lab23June3e}
\al^{m+1}\aR=\al
\left(
\sum(\nda^{v_1}\aR)\cdots(\nda^{v_j}\aR)+\BigO(Q^{n/2-2-m})
\right)
.
\end{equation}
By expanding the right-hand side of \nn{lab23June3e} using the
Leibniz rule and the formula $\al=\nda^A\nda_A$, we obtain
an expression of the form
$$\sum(\nda^{u_1}\aR)\cdots(\nda^{u_k}\aR)+\BigO(Q^{n/2-2-(m+1)})$$
plus a sum of the form
$$\sum(\al\nda^{t_0}\aR)(\nda^{t_1}\aR)\cdots(\nda^{t_s}\aR).$$
In each case, we have $u_i\leq m+1$ and $t_i\leq m$.  But by
Proposition~\ref{lab23June3a},
\[
\al\nda^{t_0}\aR
=
\sum(\nda^p\aR)(\nda^q\aR)+\BigO(Q^{n/2-3-t_0})
,
\]
where $p+q=t_0\leq m$.  Thus
\[
\al\nda^{t_0}\aR
=\displaystyle\sum(\nda^p\aR)(\nda^q\aR)
+\BigO(Q^{n/2-2-(m+1)})
.
\]
The use of the inductive assumption and Proposition~\ref{lab23June3a}
never moves an index from an $\aR$ onto a $\nda$.  $\Box$

\begin{theorem}\label{curv}
  Suppose that $n$ is even and $M$ is generic.
  Let $\h$ be an ambient metric for a conformal manifold of dimension
  $n$.  Let $t\geq 0$ and $u\geq 0$ be given, and suppose that
  $t+u\leq\frac{n}{2}-3$.  Then there is a partial contraction $\cp$,
  polynomial in $\D_A$, $\aR_{ABCD}$, $\X_{A}$, $\h_{AB}$, and its
  inverse $\h^{AB}$, such that
\begin{equation}\label{Lab17July3a}
\nda^{t}\al^{u}\aR=\cp+
\BigO(Q^{n/2-2-t-u})
.
\end{equation}
Each term of $\cp$ is of degree at least 1 in $\aR_{ABCD}$.
If, in (\ref{Lab17July3a}), $\aR$ has any free indices, then in $\cp$
these indices always appear on an $\aR$.
\end{theorem}%
\noindent{\bf Proof:}  By Proposition~\ref{lab23June3c}, we may write
\begin{equation}\label{lab23June3f}
\nda^t\al^u\aR
=\sum(\nda^{v_1}\aR)\cdots(\nda^{v_j}\aR)+\BigO(Q^{n/2-2-u-t}),
\end{equation}
where $v_i\leq t+u$ for each $i$.  If the $\aR$ on the left-hand side
of \nn{lab23June3f} has any free indices, then for each term in the
sum, these indices always appear on an $\aR$; this follows from
Proposition~\ref{lab23June3c}.  To complete the proof, we show that if
$0\leq\ell\leq\frac{n}{2}-3$, then
$\nda^{\ell}\aR=\cp+\BigO(Q^{n/2-2-\ell})$.  We use induction.
Suppose that $1\leq m\leq \frac{n}{2}-3$, and suppose that
$\nda^{\ell}\aR=\cp+\BigO(Q^{n/2-2-\ell})$ whenever $0\leq\ell\leq
m-1$.  By (\ref{lab23June3g}) we have
\begin{equation}\label{Lab30Aug3a}
\D_A\nda^{m-1}\aR=
(n-2m-4)\nda_{A}\nda^{m-1}\aR-\X_A\al\nda^{m-1}\aR
.
\end{equation}
Note that $n-2m-4>0$. 
Also observe that
each $\aR$ in \nn{Lab30Aug3a} has the same indices.  {}From
\nn{Lab30Aug3a} and Proposition~\ref{lab23June3a}, we conclude that
\[
\nda^m\aR=\D_{A}\nda^{m-1}\aR+
\X_A
\left(\rule{0mm}{5mm}
\sum(\nda^{p}\aR)(\nda^{q}\aR)+
\BigO(\rule{0mm}{4mm}Q^{n/2-3-(m-1)})
\right) ,
\]
where $p+q\leq m-1$.  Also note that if the $\aR$ on the left-hand
side of this equation has any free indices, then in each term of the
right-hand side, these indices always appear on an $\aR$.  {}From our
inductive assumption we now see that
$\nda^{m}\aR=\cp+\BigO(Q^{n/2-2-m})$. 
$\Box$

\noindent{\bf Remark:}
Theorem~\ref{curv} shows that when $n\neq 4$, an ambient partial
contraction $\nda^{t}\al^{u}\aR|_\cq$ is equivalent to a conformal
invariant which is obtained by taking a partial contraction polynomial
in $D$, $W$, $X$, $h$, and its inverse $h^{-1}$.  Moreover in each
case, via the inductive steps of the proof, one obtains the explicit
formula for the invariant as a partial contraction of these
quantities. More generally, this shows that any ``Weyl invariant''
(c.f. \cite{BEGr,F}) arising from a complete (partial) contraction of
ambient tensors of the form \nn{Lab17July3a} is contained in the space
of invariants generated by complete (partial) contractions of the
expressions polynomial in the tractor operators and fields $D$, $W$,
$X$, $h$, and $h^{-1}$. Furthermore, there is an explicit algorithm
for finding the tractor formula, given the formula for the ambient
invariant. This is a slight generalisation of a result along these
lines obtained in \cite{CapGoamb}.  $\quad\endrk$

The next proposition is a simple generalisation of results in
\cite{BrGodeRham,GJMS}.
\begin{proposition}\label{powerslaptan}
For every integer $m\geq 1$ and every ambient 
homogeneous tensor space $\act^{\Phi}(m-n/2)$,
$$
\ala^m:\act^{\Phi}(m-n/2)\to\act^\Phi(-m-n/2)
$$ 
is tangential and so determines a conformally invariant operator 
$$
\Box_m^{\alpha}: \ct^{\Phi}[m-n/2]\to\ct^\Phi[-m-n/2] .
$$
\end{proposition}
\noindent{\bf Proof:} By construction, the operators $\ala$ preserve
tensor type (tensor type with respect to pointwise ${\rm SO}(\h)$
tensor decompositions) and lower homogeneity weight by $2$. Hence $\ala^m$
maps $\act^{\Phi}(m-n/2)$ to $\act^\Phi(-m-n/2)$.

To show that $\ala^m$ acts tangentially, we calculate $\ala^mQA$ for
$A$ of homogeneity $m-2-n/2$.  Without any homogeneity assumption, we
have
\begin{equation}\label{mth}
[\ala^m,Q]=\sum_{p=0}^{m-1}\ala^{m-1-p}[\ala,Q]\ala^p.
\end{equation}
If we let \nn{mth} act on $\act^\Phi(w)$, then by \nn{abracketQ}, the
$p^{\underline{{\rm th}}}$ term on the right acts as
$2[2(w-2p)+n+2]\ala^{m-1}$.  Hence $[\ala^m,Q]$ acts as
$2m(2w-2m+n+4)\ala^{m-1}$.  This vanishes identically if $w=m-2-n/2$.
Thus $\ala^m$ is tangential on $\act^\Phi(m-n/2)$ as desired. \quad
$\Box$

The remainder of this section is concerned with obtaining tractor
formulae for the operators in the previous theorem. A key idea is to
assume that the ambient tensor field being acted on is suitably
``harmonic'' as in the following lemma. Since tangential operators do
not depend on how the field is extended off $\cq$, this involves no
loss of generality.
%
%
\begin{lemma}
\label{Lab13July3a}
Suppose $k\geq 2$ is an integer.  In the generic $n$-even case, suppose
$k\leq\frac{n}{2}-1$ or that $\alpha =0$ and $k\leq\frac{n}{2}$.  Let
$S\in\act^{\Phi}(k-\frac{n}{2})$ be given, and suppose $\ala S$ is
$\BigO(Q^{k-1})$.  Finally, let $v$, $0\leq v\leq k-1$, be given.  Then
there is a linear differential operator $\cp$ of order at most $2v$
given by a partial contraction formula polynomial in $\X_{A}$, $\D_A$,
$\aR_{ABCD}$, $\h_{AB}$, and $\h^{AB}$, such that
\begin{equation}\label{Lab13July3b}
\nda^{v}S=\cp S +\BigO(Q^{k-v}).
\end{equation}
If, on the left-hand side of (\ref{Lab13July3b}), $S$ has any free
indices, then in every term of $\cp S$, each of them appears either on
$S$ in its natural position or on $\aR$.
\end{lemma}
\noindent{\bf Proof:} We will assume that $n$ is even and $M$ is
generic.  For $v=1$, observe that by \nn{lab23June3g} and \nn{aLapl}
we have have $2(k-1)\nda S=\D S-\alpha\X\aR\#\# S+ \X\ala S $.  This
is in the required form, since $\ala S=\BigO(Q^{k-1})$.

We now proceed by induction on $v$.  Suppose
that $1\leq m<k-1$ and that (\ref{Lab13July3b}) holds for $1\leq v\leq
m$.  By (\ref{lab23June3g}) it follows that
\begin{equation}\label{Lab30Aug3c}
2(k-m-1)\nda^{m+1}S=\D\nda^m S - \alpha\X\aR\#\#\nda^m S +\X\ala\nda^{m}S.
\end{equation}
If, on the left-hand side of \nn{Lab30Aug3c}, $S$ has any free
indices, then in every term of the right-hand side of \nn{Lab30Aug3c},
each of these indices appears on an $S$ in its natural position or on
an $\aR$.  {}From the inductive assumption and the properties of $\D$,
it follows that $\D \nda^m S - \alpha\X \aR\#\#\nda^m S$ is of the form $\cp
S+\BigO(Q^{k-(m+1)})$,
where $\cp$ is as described in the
statement of the lemma.  On the other hand, by (\ref{abracketnda}),
\begin{equation}\label{Lab30Aug3d}
\begin{array}{l}
\ala\nda^{m}S=\nda^m\ala S+\sum(\nda^p\aR)(\nda^q S)
+\alpha \sum(\nda^{p+1}\aR)(\nda^{q-1} S)
\\
\hspace{20mm}
+\BigO(Q^{n/2-2-(m-1)}),
\end{array}
\end{equation}
where $p+q=m$, $p\geq 0$, and $q\geq 1$.  When we use \nn{abracketnda}
to construct \nn{Lab30Aug3d}, each index attached to $S$ on the
left-hand side of \nn{Lab30Aug3d} either remains fixed or moves onto
an $\aR$.  Note that $\nda^m\ala S$ is $\BigO(Q^{k-(m+1)})$ and that
$n/2-2-(m-1)\geq k-(m+1)$.
Thus $\ala\nda^m
S=\sum(\nda^x\aR)(\nda^y S)+\BigO(Q^{k-(m+1)})$.  Here $x+y=m$, $x\leq
m$, and $y\leq m$.  If $\alpha =0$, then we have $1\leq y$ and $x\leq
m-1$.  By Theorem~\ref{curv} and by our inductive assumption, it
follows that
\[
\ala\nda^m S=\cp S+\BigO(Q^{k-(m+1)}),
\]
where $\cp$ is as in the statement of the lemma.  $\Box$

The usefulness of Lemma~\ref{Lab13July3a} results from the next
proposition, which generalises to ambient tensors and
$\ala$-Laplacians a result of \cite{GJMS}.
%
%
\begin{proposition}
\label{Lab13July3g}
Let $k\geq 1$ be an integer.  Then for any
$T\in\act^{\Phi}(k-\frac{n}{2})$, there is an
$S\in\act^{\Phi}(k-\frac{n}{2})$ such that $T-S$ is $\BigO(Q)$ and
$\ala S$ is $\BigO(Q^{k-1})$.
\end{proposition}
\noindent{\bf Proof:} Let $w:=k-n/2$. Suppose that
$S_{m-1}\in\act^{\Phi}(w)$ is such that $T-S_{m-1}$ is $\BigO(Q)$ and
$\ala S_{m-1}=Q^{m-1} E$. (Then $E\in \act^{\Phi}(w-2m) $.) If $A\in
\act^{\Phi}(w-2m) $, then $S_m:=S_{m-1}+Q^m A \in \act^{\Phi}(w)$ and
$T-S_{m}$ is $\BigO(Q)$.  We have
$$
\ala S_m =Q^{m-1} E + \ala Q^m A.
$$  
Now
$$
\ala Q^m A = \sum_{i=0}^{m-1} Q^i[\ala,Q]Q^{m-i-1} A +\BigO(Q^{m}),
$$
and from \nn{abracketQ} and the homogeneity of $A$ and $Q$ this becomes
$$
\begin{array}{rll}
{\displaystyle \ala Q^m A}
&=& \displaystyle
 \sum_{i=0}^{m-1} 2(n+2w-4i-2) Q^{m-1} A+\BigO(Q^{m})
\vspace{1mm}
\\
&=& \displaystyle 4m(w+n/2-m)Q^{m-1}A +\BigO(Q^{m}).
\end{array}
$$
Thus if $m\neq w+n/2$
(i.e. $m\neq k$), then setting $A= -[4m(w+n/2-m)]^{-1}E$ gives $\ala
S_m=\BigO(Q^m)$.  \quad $\Box$
\\
\noindent Note that the proof establishes much
more than we require in the proposition. In particular, it shows that
the $\ala$-harmonic extension of $T|_\cq$ only fails at $\BigO(Q^k)$
and that past this the extension continues. Also, if we allow $w$ such
that $w+n/2\notin \{1,2,\cdots \}$, then for any $T \in\act^{\Phi}(w)
$ and any integer $\ell\geq 0$, there is $S\in\act^{\Phi}(w)$ such
that $T-S$ is $\BigO(Q)$ and $\ala S$ is $\BigO(Q^\ell)$.

\smallskip

\noindent{\bf Remark:} Recall that one of our central aims (at least
for $n\geq 6$) is to understand the result of applying $\afl^{n/2-2}$
to the ambient curvature $\aR$. Note that for this it would appear
that we do not need Proposition~\ref{Lab13July3g}, since by
\nn{aflaR}, the ambient curvature already has the property we require
of $S$, viz.\ that $\al_{1/2}\aR=\afl\aR =O(Q^{n/2-3})$. On the other
hand, we prefer here to treat $\afl^{n/2-2}\aR$ in two steps. First,
we derive a tractor formula for the conformally invariant operator
$\fb_{n/2-2}$ on $\ct^{2,2}[-2]$.  For this we will use
Proposition~\ref{Lab13July3g}.  This operator arises from
$\afl^{n/2-2}$ on $\act^{2,2}(-2)$.  Then finally we may apply the
operator $\fb_{n/2-2}$ to the tractor field $W$.  (See
\nn{Wform}). Proceeding in this way, we can be sure that the tractor
formula that we obtain for the ambient quantity $\afl^{n/2-2}\aR|_\cq$
is precisely the tractor formula for $\fb_{n/2-2}$ on $\ct^{2,2}[-2]$
applied to $W$. \quad \endrk

\smallskip

Next, we need to understand how powers of the $\ala$-Laplacian are
related to iterations of $\D$. We begin with a lemma which indicates
the impact of moving Laplacians to the right of $\nda$'s.
\begin{lemma}
\label{arranging}
Suppose that $n$ is even and $M$ is generic.  Let $\alpha\in\bR$,
$w\in\bR$, and $T\in\act^{\Phi}(w)$ be given.  Let
\begin{equation}\label{SDefinition}
S=\ala^{t_1}\nda^{u_1}\cdots\ala^{t_p}\nda^{u_p}T,
\end{equation}
where $t_i+u_i\geq 1$ for each $i$.  Suppose that
$k:=\sum_{i=1}^p(t_i+u_i)\leq \frac{n}{2}-1$.  Then
\begin{equation}\label{SFormula}
S=\sum
(\nda^{v_1}\ala^{w_1}\aR)\cdots(\nda^{v_q}\ala^{w_q}\aR)
(\nda^{v_{q+1}}\ala^{w_{q+1}}T)
+\BigO(Q^{n/2-k})
,
\end{equation}
where $v_j+w_j\leq k$ for each $j$. 
If $T$ has any free indices in \nn{SDefinition}, then
in \nn{SFormula} these indices appear either on $T$ in their original
position or on an $\aR$.
\end{lemma}

\noindent{\bf Proof:} We proceed by induction on $k$.  Suppose that
$1\leq m\leq \frac{n}{2}-2$.  Suppose the result holds whenever $1\leq
k\leq m$, and let $S$ be as in \nn{SDefinition} with $k=m+1$.  If
$t_1=0$, then by our inductive assumption we see immediately that
\nn{SFormula} holds modulo $\BigO(Q^{n/2-(m+1)})$.  On the other hand,
suppose $t_1>0$.  Then by our inductive assumption,
\[
S=
\ala\left(\sum(\nda^{v_1}\ala^{w_1}\aR)\cdots
(\nda^{v_q}\ala^{w_q}\aR)(\nda^{v_{q+1}}\ala^{w_{q+1}}T)
+\BigO(Q^{n/2-m})\right),
\]
where $v_j+w_j\leq m$ for each $j$.  Suppose we use the Leibniz rule
to expand
\[
\al(\nda^{v_1}\ala^{w_1}\aR)\cdots
(\nda^{v_q}\ala^{w_q}\aR)(\nda^{v_{q+1}}\ala^{w_{q+1}}T).
\]
Then each term in the resulting sum will contain two factors of the
form $\nda^{v_j+1}\ala^{w_j}\mbox{\boldmath{$P$}}$ or one factor of
the form $\al\nda^{v_j}\ala^{w_j}\mbox{\boldmath{$P$}}$, where
$\mbox{\boldmath{$P$}}$ denotes $\aR$ or $T$ in each case.  But
\[
\al\nda^{v_j}\ala^{w_j}\mbox{\boldmath{$P$}}
  =
\ala\nda^{v_j}\ala^{w_j}\mbox{\boldmath{$P$}}
  -
\alpha\aR\#\#\nda^{v_j}\ala^{w_j}\mbox{\boldmath{$P$}},
\]
and by \nn{abracketnda} we may write $\ala\nda^{v_j}\ala^{w_j}
\mbox{\boldmath{$P$}}$ in the form
\[
\nda^{v_j}\ala^{w_j+1}\mbox{\boldmath{$P$}}+
\sum(\nda^{v_{\ell}^{'}}\aR)\nda^{v_{\ell}^{''}}\ala^{w_j}
\mbox{\boldmath{$P$}}+
\BigO(Q^{n/2-(m+1)})
.
\]
Here $v_{\ell}^{'}+v_{\ell}^{''}=v_j$.  When we use \nn{abracketnda},
any given index attached to $\mbox{\boldmath{$P$}}$ either remains
fixed or moves onto an $\aR$.  This completes the induction.  $\Box$
%
%
\begin{lemma}
\label{Lab11July3d}
Suppose $\ell$ is an integer and $\ell\geq 1$.  In the generic
$n$-even case, suppose also that $\ell\leq\frac{n}{2}-1$.  Let
$T\in\act^{\Phi}(\ell-\frac{n}{2})$ be given.  Then
\begin{equation}\label{Lab3July3a}
\begin{array}{lll}
\lefteqn{\ala^{\ell-1}\D T =}
&&
\\
&& -\X\ala^{\ell}T+
\sum(\nda^{v_1}\al^{w_1}\aR)\cdots
(\nda^{v_p}\al^{w_p}\aR)
(\nda^{v_{p+1}}\ala^{w_{p+1}}T)
\\
&& +\alpha \X \sum(\nda^{r_1}\al^{s_1}\aR)\cdots
(\nda^{r_q}\al^{s_q}\aR)
(\nda^{r_{q+1}}\ala^{s_{q+1}}T) +\BigO(Q).
\end{array}
\end{equation}
Here $v_i+w_i\leq \ell-1$ for $1\leq i\leq p+1$, and $r_i+s_i\leq
\ell-1$ for $1\leq i\leq q+1$.  If $\alpha=0$, then $v_i+w_i\leq
\ell-2$ for $1\leq i\leq p$, and $v_{p+1}+w_{p+1}\leq\ell-1$.  If, on
the left-hand side of (\ref{Lab3July3a}), $T$ has any free indices,
then on the right-hand side these indices always appear on $\aR$ or in
their natural positions on $T$.
\end{lemma}
\noindent{\bf Proof:} Suppose that $n$ is even and $M$ is generic.  If
$\ell=1$, the result follows from \nn{lab23June3g}.  Now suppose
that $\ell\geq 2$.  {}From (\ref{lab23June3g}) and \nn{abracketalX} we
have

\[
\begin{array}{lll}
\lefteqn{\displaystyle\ala^{\ell-1}\D_{A}T}
&&
\vspace{2mm} \\ &&
=2(\ell-1)\ala^{\ell-1}\nda_AT-\ala^{\ell-1}\X_A\ala T+\alpha
\ala^{\ell-1}\X_A \aR\#\# T
\vspace{2mm}
\\
&&
=\displaystyle
2(\ell-1)\ala^{\ell-1}\nda_A T
-[\ala^{\ell-1},\X_A]\ala T
-\X_A\ala^{\ell-1}\ala  T
\vspace{2mm}
\\
&&
\hspace*{4mm} +\alpha  [\ala^{\ell-1},\X_A] \aR\#\# T 
+\alpha  \X_A \ala^{\ell-1}\aR\#\# T 
\vspace{2mm} \\
&&
=\displaystyle
-\X_A\ala^{\ell}T
+2(\ell-1)\ala^{\ell-1}\nda_A T
-
\rule{0mm}{5mm}
\sum_{i=0}^{\ell-2}\ala^{\ell-2-i}[\ala,\X_A]\ala^{i}
\rule{0mm}{5mm}
\ala T
\vspace{2mm}
\\
&&
\displaystyle\hspace*{4mm} +\alpha      \left(
\rule{0mm}{5mm}
\right.
\sum_{i=0}^{\ell-2}\ala^{\ell-2-i}[\ala,\X_A]\ala^{i}
\left.
\rule{0mm}{5mm}
\right)       \aR\#\# T 
+\alpha  \X_A \ala^{\ell-1}\aR\#\# T 
\vspace{2mm}
\\
&&
=\displaystyle
-\X_A\ala^{\ell}T
+2(\ell-1)\ala^{\ell-1}\nda_A T
-
2
\rule{0mm}{5mm}
\sum_{i=0}^{\ell-2}\ala^{\ell-2-i}\nda_A\ala^{i}
\rule{0mm}{5mm}
\ala T
\vspace{2mm}
 \\
&&
\displaystyle\hspace*{4mm} +2\alpha      \left(
\rule{0mm}{5mm}
\right.
\sum_{i=0}^{\ell-2}\ala^{\ell-2-i}\nda_A\ala^{i}
\left.
\rule{0mm}{5mm}
\right)       \aR\#\# T 
+\alpha  \X_A \ala^{\ell-1}\aR\#\# T .
\end{array}
\]
Each of the original indices on $T$ remains fixed in the above
calculation except in the terms of $\aR\#\# T$, where it may either
remain in its original position on $T$ or move onto an $\aR$.  By
(\ref{abracketnda}), we may re-express this in the form
\begin{equation}\label{Lab30Aug3e}
\begin{array}{lll}
\lefteqn{\displaystyle\ala^{\ell-1}\D T=}
&&
\\
&&
-\X\ala^\ell T
+\sum\ala^{s_j}\aR\nda \ala^{t_j}T+\alpha \sum\ala^{s_j}(\nda \aR)\ala^{t_j}T
\vspace{2mm} \\ 
&& 
+\alpha  \ala^{\ell-2}\nda \aR\#\# T + 
\alpha \sum\ala^{p_i}\aR\nda \ala^{q_i}\aR\#\#  T
\vspace{2mm}
\\
&&
 +\alpha \sum\ala^{p_i}(\nda\aR)\ala^{q_i}\aR\#\#  T
 + \alpha \X \ala^{\ell-1}\aR\#\# T+\BigO(Q),
\end{array}
\end{equation}
where $s_j +t_j=\ell-2$ for each $j$ and $p_i+q_i=\ell-3$ for each
$i$.  When we use \nn{abracketnda} to construct \nn{Lab30Aug3e}, each
index on $T$ or $\aR$ either remains fixed or moves onto an $\aR$. In
the right-hand side of \nn{Lab30Aug3e} the coefficient of $\X\ala^\ell
T$ is exact.  Otherwise, no attempt has been made to present the
coefficients precisely. At this point we need only the general form of
the expression. Where there is a coefficient $\alpha$ presented, this
means, as usual, that all terms of this form appear with coefficient a
multiple of $\alpha$.
 
For ambient tensors $U$ and $V$,
\[
\ala UV=
(\al U)V+(\nda U)\nda V+U\ala V +\aR UV
.
\]
Thus by using the definition of $\ala$ together with the Leibniz rule,
we may re-express the right-hand side of \nn{Lab30Aug3e} in the form
given on the right-hand side of \nn{Lab3July3a}, except that on each
$\aR$ or $T$, the operators $\nda$, $\al$, and $\ala$ may not be in
the order given in \nn{Lab3July3a}.  But by Lemma~\ref{arranging}, we
may indeed re-express the right-hand side of \nn{Lab30Aug3e} in the
form given on the right-hand side of \nn{Lab3July3a}.  In doing this,
we may move an index that was originally attached to an $\aR$ or a
$T$, but we always move the index onto an $\aR$.  In the new
expression, we have $v_i+w_i\leq \ell-1$ for $1\leq i\leq p+1$ and
$r_i+s_i\leq \ell-1$ for $1\leq i\leq q+1$; this follows from
Lemma~\ref{arranging}.  In the $\alpha=0$ case, the fact that
$v_i+w_i\leq \ell-2$ for $1\leq i\leq p$ follows from the fact that
\nn{Lab30Aug3e} simplifies to
$$
\al^{\ell-1}\D T =-\X\al^\ell T
+\sum\al^{s_j}\aR\nda \al^{t_j}T + \BigO(Q)
$$ when $\alpha=0$.  $\Box$

We are now ready to show that the powers of the $\ala$-Laplacian can
be re-expressed as a sum of compositions of tangential operators. 
%
\begin{proposition}
\label{Lab11July3a}
Suppose $k\geq1$ is an integer.  Let $w=k-\frac{n}{2}$, and let
$V\in\act^{\Phi}(w)$ be given.  In the generic $n$-even case, suppose
that $k\leq \frac{n}{2}-2$, or $\alpha =0$ and $k\leq \frac{n}{2}-1$,
or $\act^\Phi(w)=\act(w)$ and $k\leq
\frac{n}{2}-1$, or $\act^\Phi(w)=\act^0(w)$
and $k\leq n/2$.  Then
\begin{equation}\label{Lab3July3b}
(-1)^{k-1}\X_{A_1}\cdots\X_{A_{k-1}}\ala^{k}V=
\al\D_{A_1}\cdots\D_{A_{k-1}}V+\cp V
+
\BigO(Q)
,
\end{equation}
where $\cp$ is a linear differential operator of order less than $2k$
given as a partial contraction polynomial in $\X_{A}$, $\D_A$,
$\aR_{ABCD}$, $\h_{AB}$, and $\h^{AB}$.  If, in (\ref{Lab3July3b}),
$V$ has any free indices, then for every term of $\cp V$, these
indices appear either on $\aR$ or in their natural position on $V$.
The indices $A_i$ are \textit{not} skew-symmetrised.
\end{proposition}
%

\noindent{\bf Proof:} The case of $V\in \act^0(w)$ is treated in
\cite{GoPet}.  For the remaining cases, we assume, as usual, that we
are in the generic $n$-even setting.

We begin with the case $k\leq\frac{n}{2}-2$ and the case $\alpha=0$,
$k\leq\frac{n}{2}-1$, and we use induction on $k$.  Suppose that
$1\leq m\leq\frac{n}{2}-3$ or that $\alpha=0$ and $1\leq
m\leq\frac{n}{2}-2$, and suppose the result holds whenever $k=m$.  Let
$V\in\act^{\Phi}(m+1-\frac{n}{2})$.  By Proposition~\ref{Lab13July3g},
there exists an $S\in\act^{\Phi}(m+1-\frac{n}{2})$ such that $V-S$ is
$\BigO(Q)$ and $\ala S$ is $\BigO(Q^m)$.  Then by our inductive
assumption,
\begin{equation}\label{Lab11July3b}
\begin{array}{lll}
\lefteqn{(-1)^{m-1}\X_{A_1}\cdots\X_{A_{m-1}}\ala^{m}(\D_{A_m}S)=}
&&\vspace{2mm}
\\
&&
\al\D_{A_1}\cdots\D_{A_{m-1}}(\D_{A_m}S)
+\cp S+\BigO(Q)
,
\end{array}
\end{equation}
where $\cp$ is of order less than $2m$.
If, on the left-hand side of \nn{Lab11July3b}, $S$ has any free
indices, then in each term of $\cp S$, these indices appear either on
$\aR$ or in their natural position on $S$.  Now apply
Lemma~\ref{Lab11July3d} with $\ell=m+1$ and $T=S$.  We conclude that
\begin{equation}\label{Lab11July3e}
\begin{array}{lll}
\lefteqn{\ala^{m}\D_{A_m}S=}&&\vspace{2mm}
\\
&&
-\X_{A_m}\ala^{m+1}S+\BigO(Q)
\vspace{2mm}
\\
&&
+
\sum
(\nda^{v_1}\al^{w_1}\aR)\cdots
(\nda^{v_p}\al^{w_p}\aR)
(\nda^{v_{p+1}}\ala^{w_{p+1}}S)
\vspace{2mm}
\\
&&
+
\alpha\X\sum(\nda^{r_1}\al^{t_1}\aR)\cdots
(\nda^{r_q}\al^{t_q}\aR)(\nda^{r_{q+1}}\ala^{t_{q+1}}S)
.
\end{array}
\end{equation}
Here $v_i+w_i\leq m$ for $1\leq i\leq p+1$, and $r_i+t_i\leq m$ for
$1\leq i\leq q+1$.  If $\alpha=0$, then $v_i+w_i\leq m-1$ for $1\leq
i\leq p$ and $v_{p+1}+w_{p+1}\leq m$.  If, on the left-hand side of
\nn{Lab11July3e}, $S$ has any free indices, then on the right-hand
side of this equation these indices appear on $\aR$ or in their
natural positions on $S$.  Since $\ala S$ is $\BigO(Q^m)$, we may
assume that $w_{p+1}=t_{q+1}=0$ in (\ref{Lab11July3e}).  Thus by
Theorem~\ref{curv} and Lemma~\ref{Lab13July3a}, we have
\begin{equation}\label{Lab11July3f}
\ala^{m}\D_{A_m}S=-\X_{A_m}\ala^{m+1}S+\cp S+\BigO(Q)
.
\end{equation}%
Since $v_{p+1}\leq m$ and $r_{q+1}\leq m$ in \nn{Lab11July3e}, it follows that the order of $\cp$ is at most $2m$ in \nn{Lab11July3f}.
If, in \nn{Lab11July3f}, $S$ has free indices, then in $\cp S$ these
appear either on $\aR$ or in their natural positions on $S$.  {}From
(\ref{Lab11July3b}) and (\ref{Lab11July3f}) it now follows that
\begin{equation}\label{Lab11July3g}
(-1)^m \X_{A_1}\cdots\X_{A_m}\ala^{m+1}S
=
\al\D_{A_1}\cdots\D_{A_m}S
+\cp S
+\BigO(Q)
.
\end{equation}
But $\D_A$ acts tangentially along $\cq$, and $\al$ acts tangentially
on fields homogeneous of degree $1-n/2$.  Thus
$\al\D_{A_1}\cdots\D_{A_m} +\cp$ acts tangentially on $S$.  By
Proposition~\ref{powerslaptan}, $\ala^{m+1}$ also acts tangentially on
$S$, and so we may replace $S$ with $V$ on both sides of
(\ref{Lab11July3g}).  This completes the induction.

Finally, suppose that $\act^{\Phi}(w)=\act(w)$.  By the Ricci
flatness of the ambient metric, it follows that $\aR\#\#V$ is
$\BigO(Q^{n/2-1})$.  Thus for $1\leq k\leq\frac{n}{2}-1$ we see that
$\ala^kV=\al^k V+\BigO(Q)$, and the result follows from the case
$\alpha=0$.  $\Box$

We are now ready to prove Theorem~\ref{powerslap} and at the same
time describe tractor formulae for the operators $\Box^\alpha_m$.  We
begin with the tractor formulae.
\begin{theorem}\label{powersnatural}
Via the algorithm implicit in the inductive steps above, the operators
$\Box_m^{\alpha}$ have tractor formulae (for $m$ in the ranges given
in Theorem~\ref{powerslap}) as follows:
%
\begin{equation}\label{powerform}
\begin{array}{lll}
\lefteqn{\displaystyle(-1)^{m-1}X_{A_1}\cdots
X_{A_{m-1}} \Box^\alpha_m U }
&&\vspace{2mm}
\\
&&
\displaystyle
=
\Box D_{A_1}\cdots D_{A_{m-1}}
U
+\cp^{\Phi,m}_{A_1\cdots A_{m-1}} U
,
\end{array}
\end{equation}
%
where the differential operator $\cp^{\Phi,m}$ is a partial
contraction polynomial in $X$, $D$, $W$, $h$, and $h^{-1}$.  Thus for
$m\neq n/2$,
\begin{equation}\label{powerformtoo}
\begin{array}{lll}
\lefteqn{\displaystyle(m-1)!\left(\Pi_{i=2}^{m}(n-2i)\right)
\Box^\alpha_m U =}
&&\vspace{2mm}
\\
&&
\displaystyle
D^{A_{m-1}}\cdots D^{A_1}
\Box D_{A_1}\cdots D_{A_{m-1}}U
\vspace{2mm}
\\
&&
\displaystyle
+
D^{A_{m-1}}\cdots D^{A_1} 
\cp^{\Phi,m}_{A_1\cdots A_{m-1}} U
.
\end{array}
\end{equation}
The indices attached to $U$ on the left-hand side appear, in each term
of $\cp^{\Phi,m} U$, on $U$ in their original position or on $W$.  The
indices $A_i$ in \nn{powerform} and \nn{powerformtoo} are \textit{not}
skew-symmetrised.
\end{theorem}
\noindent{\bf Proof of theorems~\ref{powerslap} and
 \ref{powersnatural}:} Recall that $\D$ descends to $D$, and
 $\al:\act^{\Phi}(1-n/2)\rightarrow \act^{\Phi}(-1-n/2)$ descends to
 the generalised conformal Laplacian operator
 $\Box$.  (See
 \nn{clapn}.)  Thus \nn{powerform} is an immediate consequence of
 Proposition~\ref{Lab11July3a}.  {}From this the claims of naturality
 are immediate from the naturality of $X$, $\Box$, $D$, $W$, $h$, and
 $h^{-1}$. That the $\Box^\alpha_m$ have leading term $\Delta^m$
 follows easily from the expression \nn{Dform} for $D$ and the
 identities \nn{connids} for the tractor connection.  Then note that
 \nn{powerformtoo} follows from
 (\ref{powerform}) and (\ref{Lab13July3d}).
\quad $\Box$

\subsection{Calculating explicit formulae; examples} \label{examples}

One can easily compute explicit formulae for the obstruction tensors
in low dimensions.  {}From the proof of theorems~\ref{main2} and
\ref{main1}, above, we know that in dimension 4, $\cB_{ab}$ is simply
$-\frac{1}{2}B_{ab}$, where $B_{ab}$ is the Bach tensor as given in \nn{one}.  

In dimension 6, we have $m=1$, and the relevant ambient
operator from Proposition~\ref{powerslaptan} is $\afl:
\act^{2,2}(-2)\to \act^{2,2}(-4)$, which descends to
\begin{equation}\label{6Lap}
\Box  + \frac{1}{4}W\#\# =:\fb_1: \ct^{2,2}[-2]\to \ct^{2,2}[-4] .
\end{equation}
The left-hand side of \nn{6Lap} is the tractor formula for $\fb_1$.
By Theorem \ref{main2}, applying this to $W$ (see \nn{Wform}) yields
the obstruction tensor via the identity \nn{tractor4obstrn}. That is,
$$
2^6 X_{A_1}Z_{A_2}{}^aX_{B_1}Z_{B_2}{}^b \cB^{6}_{ab} =
 \Box W  + \frac{1}{4}W\#\#W,
$$
where we have used the fact that $k(6)= 2^6$.  Thus 
\[
\begin{array}{c}
64 \cB_{ce} X_{[B} Z_{C]}{}^c Z_{[E}{}^eX_{D]}=
\vspace{2mm}
\\
\Box W_{BCDE}-
W^A{}_{CB}{}^{F}W_{FADE}-
W^{A}{}_{CD}{}^{F}W_{BAFE}-
W^{A}{}_{CE}{}^{F}W_{BADF}

.
\end{array}
\]
But $4Y^B Y^D Z^{C}{}_{a}Z^{E}{}_{b} X_{[B}B_{C][E}X_{D]}= \cB_{ab} $.
Thus in any conformal scale, $\cB_{ab}$ is given by the following
formula:
\begin{equation}\label{BSix}
\begin{array}{lll}
\lefteqn{{\textstyle\frac{1}{16}}
Y^B Y^D Z^C{}_{a}Z^E{}_{b}\left(\Box W_{BCDE}
\rule{0mm}{4mm}
\right.}
&&\vspace{3mm}
\\
&&
\left.
\rule{0mm}{4mm}
-W^A{}_{CB}{}^F W_{FADE}
-W^A{}_{CD}{}^F W_{BAFE}
-W^A{}_{CE}{}^F W_{BADF}
\right)
.
\end{array}
\end{equation}
If one expands using \nn{connids}, \nn{Wform}, and the definitions of
$\Box$ and the tractor metric, it is an entirely mechanical process to
rewrite \nn{BSix} in terms of the Levi-Civita connection and its
curvature (with metric contractions).  A computation using this
process together with \textit{Mathematica} and J.~Lee's \textit{Ricci}
software package (\cite{Lee}) shows that
\begin{equation}\label{B6}
\begin{array}{lll}
\lefteqn{\textstyle\cB^{6}_{ab}=
\frac{1}{16}\Delta B_{ab}
-\frac{1}{4}\J B_{ab}
+\frac{1}{8}B_{cd}C_{a}{}^{c}{}_{b}{}^{d}
-
\frac{1}{2}\V_{cd}\nd^{c}\CY_{(ab)}{}^{d}}
\vspace{2mm}
\\
&&
+\frac{1}{4}\CY_{cad}\CY^{c}{}_{b}{}^{d}
-\frac{1}{2}\CY_{cad}\CY^{d}{}_b{}^c
-\frac{1}{4}\CY_{(ab)c}\nd^c\J
+\frac{1}{4}\V_{cd}\V^{d}{}_e C_a{}^c{}_b{}^e
,
\end{array}
\end{equation}
%
%
where $\CY$ and $B$ are respectively the Cotton and Bach tensors
as given in \nn{cot} and \nn{Bachform}.  This formula for
$\cB^{6}_{ab}$ agrees up to a constant factor with the formula given
by Graham and Hirachi in \cite{GrH}.

In dimension 8, we find that $\cB_{ab}=\frac{1}{384}\SymTen_{(ab)}$,
where $\SymTen_{ab}$ is as given in Figure~\ref{B8}.
\begin{figure}
\[
\begin{array}{l}
-
\tensor{\Bach}{\down{a}\down{b}\down{|}\down{c}\up{c}\down{d}\up{d}}
%
%
+10 
\tensor{\Bach}{\down{a}\down{b}\down{|}\down{c}\up{c}} \J
%
%
-28 
\tensor{\Bach}{\down{a}\down{b}\down{|}\down{c}\down{d}} 
\tensor{\V}{\up{c}\up{d}}
%
%
+24 
\tensor{\Bach}{\down{a}\down{c}\down{|}\down{b}\down{d}} 
\tensor{\V}{\up{c}\up{d}}
\\
%
%
-4 
\tensor{\Bach}{\down{c}\down{d}\down{|}\down{e}\up{e}} 
\tensor{\Weyl}{\down{a}\up{c}\down{b}\up{d}}
%
%
-24 
\tensor{\Bach}{\down{a}\down{c}\down{|}\down{d}} 
\tensor{\V}{\down{b}\up{c}\down{|}\up{d}}
%
%
-24 
\tensor{\Bach}{\down{c}\down{d}\down{|}\down{a}} 
\tensor{\V}{\down{b}\up{c}\down{|}\up{d}}
%
%
+56 
\tensor{\Bach}{\down{a}\down{c}\down{|}\down{d}} 
\tensor{\V}{\down{b}\up{d}\down{|}\up{c}}
\\
%
%
-6 
\tensor{\Bach}{\down{a}\down{b}\down{|}\down{c}} 
\tensor{\J}{\down{|}\up{c}}
%
%
+12 
\tensor{\Bach}{\down{a}\down{c}\down{|}\down{b}} 
\tensor{\J}{\down{|}\up{c}}
%
%
+24 
\tensor{\Bach}{\down{c}\down{d}\down{|}\down{a}} 
\tensor{\V}{\up{c}\up{d}\down{|}\down{b}}
%
%
-32 
\tensor{\Bach}{\down{a}\down{c}\down{|}\down{d}} 
\tensor{\V}{\up{c}\up{d}\down{|}\down{b}}
\\
%
%
-4 
\tensor{\Bach}{\down{c}\down{d}\down{|}\down{e}} 
\tensor{\Weyl}{\down{a}\up{c}\down{b}\up{d}\down{|}\up{e}}
%
%
+4 
\tensor{\Bach}{\down{a}\down{b}} 
\tensor{\J}{\down{|}\down{c}\up{c}}
%
%
-16 
\tensor{\Bach}{\down{c}\down{d}} 
\tensor{\V}{\up{c}\up{d}\down{|}\down{a}\down{b}}
%
%
-40 
\tensor{\Bach}{\down{c}\down{d}} 
\tensor{\V}{\down{a}\down{b}\down{|}\up{c}\up{d}}
\\
%
%
+56 
\tensor{\Bach}{\down{c}\down{d}} 
\tensor{\V}{\down{a}\up{c}\down{|}\down{b}\up{d}}
%
%
-8 
\tensor{\Bach}{\down{a}\down{c}} 
\tensor{\Bach}{\down{b}\up{c}}
%
%
+3 
\tensor{\Bach}{\down{c}\down{d}} 
\tensor{\Bach}{\up{c}\up{d}} 
\tensor{\bg}{\down{a}\down{b}}
%
%
-24 
\tensor{\Bach}{\down{a}\down{b}} \J^{2}
%
%
-64 
\tensor{\Bach}{\down{a}\down{c}} 
\tensor{\V}{\down{b}\down{d}} 
\tensor{\V}{\up{c}\up{d}}
\\
%
%
+76 
\tensor{\Bach}{\down{a}\down{b}} 
\tensor{\V}{\down{c}\down{d}} 
\tensor{\V}{\up{c}\up{d}}
%
%
+28 
\tensor{\Bach}{\down{c}\down{d}} 
\tensor{\bg}{\down{a}\down{b}} 
\tensor{\V}{\up{c}\down{e}} 
\tensor{\V}{\up{d}\up{e}}
%
%
+16 
\tensor{\Bach}{\down{c}\down{d}} \J 
\tensor{\Weyl}{\down{a}\up{c}\down{b}\up{d}}
%
%
+32 
\tensor{\Bach}{\down{c}\down{d}} 
\tensor{\V}{\down{a}\down{e}} 
\tensor{\Weyl}{\down{b}\up{c}\up{d}\up{e}}
\\
%
%
-24 
\tensor{\Bach}{\down{c}\down{d}} 
\tensor{\V}{\up{c}\down{e}} 
\tensor{\Weyl}{\down{a}\up{d}\down{b}\up{e}}
%
%
+4 
\tensor{\Bach}{\down{c}\down{d}} 
\tensor{\Weyl}{\down{a}\down{e}\up{c}\down{i}} 
\tensor{\Weyl}{\down{b}\up{e}\up{d}\up{i}}
%
%
-8 
\tensor{\Bach}{\down{c}\down{d}} 
\tensor{\Weyl}{\down{a}\down{e}\up{c}\down{i}} 
\tensor{\Weyl}{\down{b}\up{i}\up{d}\up{e}}
%
%
-8 
\tensor{\CY}{\down{a}\down{c}\down{b}} 
\tensor{\J}{\down{|}\down{d}\up{d}\up{c}}
\\
%
%
-32 
\tensor{\CY}{\down{a}\down{c}\down{b}\down{|}\down{d}\down{e}} 
\tensor{\V}{\up{c}\up{d}\down{|}\up{e}}
%
%
-16 
\tensor{\CY}{\down{a}\down{c}\down{d}\down{|}\down{e}} 
\tensor{\CY}{\down{b}\up{c}\up{d}\down{|}\up{e}}
%
%
+16 
\tensor{\CY}{\down{c}\down{d}\down{a}\down{|}\down{e}} 
\tensor{\CY}{\up{c}\up{d}\down{b}\down{|}\up{e}}
%
%
-32 
\tensor{\CY}{\down{a}\down{c}\down{b}\down{|}\down{d}} 
\tensor{\J}{\down{|}\up{c}\up{d}}
\\
%
%
+32 
\tensor{\CY}{\down{c}\down{a}\down{d}\down{|}\down{e}} 
\tensor{\V}{\down{b}\up{e}} 
\tensor{\V}{\up{c}\up{d}}
%
%
-64 
\tensor{\CY}{\down{a}\down{b}\down{c}\down{|}\down{d}} 
\tensor{\V}{\up{c}\up{d}} \J
%
%
-128 
\tensor{\CY}{\down{a}\down{c}\down{d}\down{|}\down{e}} 
\tensor{\V}{\down{b}\up{d}} 
\tensor{\V}{\up{c}\up{e}}
%
%
-128 
\tensor{\CY}{\down{c}\down{a}\down{d}\down{|}\down{e}} 
\tensor{\V}{\down{b}\up{d}} 
\tensor{\V}{\up{c}\up{e}}
\\
%
%
-608 
\tensor{\CY}{\down{a}\down{c}\down{b}\down{|}\down{d}} 
\tensor{\V}{\up{c}\down{e}} 
\tensor{\V}{\up{d}\up{e}}
%
%
-32 
\tensor{\CY}{\down{c}\down{a}\down{d}\down{|}\down{b}} 
\tensor{\V}{\up{c}\down{e}} 
\tensor{\V}{\up{d}\up{e}}
%
%
+32 
\tensor{\CY}{\down{a}\down{c}\down{d}\down{|}\down{e}} 
\tensor{\V}{\up{e}\down{i}} 
\tensor{\Weyl}{\down{b}\up{c}\up{d}\up{i}}
\\
%
%
+32 
\tensor{\CY}{\down{c}\down{a}\down{d}\down{|}\down{e}} 
\tensor{\V}{\up{e}\down{i}} 
\tensor{\Weyl}{\down{b}\up{c}\up{d}\up{i}}
%
%
+32 
\tensor{\CY}{\down{a}\down{c}\down{d}\down{|}\down{e}} 
\tensor{\V}{\up{d}\down{i}} 
\tensor{\Weyl}{\down{b}\up{c}\up{e}\up{i}}
%
%
+32 
\tensor{\CY}{\down{c}\down{a}\down{d}\down{|}\down{e}} 
\tensor{\V}{\up{d}\down{i}} 
\tensor{\Weyl}{\down{b}\up{c}\up{e}\up{i}}
\\
%
%
-64 
\tensor{\CY}{\down{a}\down{b}\down{c}\down{|}\down{d}} 
\tensor{\V}{\down{e}\down{i}} 
\tensor{\Weyl}{\up{c}\up{e}\up{d}\up{i}}
%
%
+32 
\tensor{\V}{\down{c}\down{d}} 
\tensor{\V}{\down{e}\down{i}} 
\tensor{\Weyl}{\down{a}\up{c}\down{b}\up{e}\down{|}\up{d}\up{i}}
%
%
+32 
\tensor{\V}{\down{c}\down{d}} 
\tensor{\J}{\down{|}\down{e}} 
\tensor{\Weyl}{\down{a}\up{c}\down{b}\up{e}\down{|}\up{d}}
\\
%
%
+32 
\tensor{\CY}{\down{c}\down{d}\down{e}} 
\tensor{\V}{\up{d}\down{i}} 
\tensor{\Weyl}{\down{a}\up{c}\down{b}\up{e}\down{|}\up{i}}
%
%
+32 
\tensor{\CY}{\down{c}\down{d}\down{e}} 
\tensor{\V}{\up{d}\down{i}} 
\tensor{\Weyl}{\down{a}\up{c}\down{b}\up{i}\down{|}\up{e}}
%
%
+64 
\tensor{\CY}{\down{a}\down{c}\down{d}} 
\tensor{\V}{\down{e}\down{i}} 
\tensor{\Weyl}{\down{b}\up{c}\up{d}\up{e}\down{|}\up{i}}
\\
%
%
+64 
\tensor{\CY}{\down{c}\down{a}\down{d}} 
\tensor{\V}{\down{e}\down{i}} 
\tensor{\Weyl}{\down{b}\up{c}\up{d}\up{e}\down{|}\up{i}}
%
%
+8 
\tensor{\J}{\down{|}\down{c}} 
\tensor{\J}{\down{|}\down{d}} 
\tensor{\Weyl}{\down{a}\up{c}\down{b}\up{d}}
%
%
-16 
\tensor{\V}{\down{c}\down{d}\down{|}\down{e}} 
\tensor{\V}{\up{c}\down{i}\down{|}\up{e}} 
\tensor{\Weyl}{\down{a}\up{d}\down{b}\up{i}}
\\
%
%
+32 
\tensor{\CY}{\down{c}\down{a}\down{d}} 
\tensor{\J}{\down{|}\down{e}} 
\tensor{\Weyl}{\down{b}\up{c}\up{d}\up{e}}
%
%
-32 
\tensor{\CY}{\down{c}\down{a}\down{d}} 
\tensor{\J}{\down{|}\down{e}} 
\tensor{\Weyl}{\down{b}\up{e}\up{c}\up{d}}
%
%
-16 
\tensor{\CY}{\down{c}\down{d}\down{e}} 
\tensor{\CY}{\up{c}\up{d}\down{i}} 
\tensor{\Weyl}{\down{a}\up{i}\down{b}\up{e}}
\\
%
%
+32 
\tensor{\CY}{\down{c}\down{d}\down{e}} 
\tensor{\CY}{\up{d}\up{c}\down{i}} 
\tensor{\Weyl}{\down{a}\up{i}\down{b}\up{e}}
%
%
-32 
\tensor{\CY}{\down{a}\down{c}\down{d}} 
\tensor{\CY}{\down{e}\up{d}\down{i}} 
\tensor{\Weyl}{\down{b}\up{e}\up{c}\up{i}}
%
%
-32 
\tensor{\CY}{\down{c}\down{a}\down{d}} 
\tensor{\CY}{\up{d}\down{e}\down{i}} 
\tensor{\Weyl}{\down{b}\up{e}\up{c}\up{i}}
\\
%
%
-32 
\tensor{\CY}{\down{c}\down{a}\down{d}} 
\tensor{\CY}{\down{e}\down{i}\up{c}} 
\tensor{\Weyl}{\down{b}\up{i}\up{d}\up{e}}
%
%
%
+64 
\tensor{\CY}{\down{c}\down{a}\down{d}} 
\tensor{\CY}{\down{e}\down{b}\down{i}} 
\tensor{\Weyl}{\up{c}\up{d}\up{e}\up{i}}
%
%
-32 
\tensor{\CY}{\down{c}\down{a}\down{d}} 
\tensor{\CY}{\down{e}\down{b}\down{i}} 
\tensor{\Weyl}{\up{c}\up{e}\up{d}\up{i}}
%
%
-64 
\tensor{\CY}{\down{a}\down{c}\down{d}} 
\tensor{\V}{\down{b}\up{d}} 
\tensor{\J}{\down{|}\up{c}}
\\
%
%
-64 
\tensor{\CY}{\down{c}\down{a}\down{d}} 
\tensor{\V}{\down{b}\up{d}} 
\tensor{\J}{\down{|}\up{c}}
%
%
-32 
\tensor{\CY}{\down{a}\down{b}\down{c}} \J 
\tensor{\J}{\down{|}\up{c}}
%
%
-16 
\tensor{\CY}{\down{c}\down{d}\down{a}} 
\tensor{\V}{\up{c}\up{d}} 
\tensor{\J}{\down{|}\down{b}}
%
%
-224 
\tensor{\CY}{\down{a}\down{c}\down{b}} 
\tensor{\V}{\up{c}\down{d}} 
\tensor{\J}{\down{|}\up{d}}
\\
%
%
-96 
\tensor{\CY}{\down{c}\down{a}\down{d}} 
\tensor{\CY}{\down{e}\up{c}\up{d}} 
\tensor{\V}{\down{b}\up{e}}
%
%
-192 
\tensor{\CY}{\down{c}\down{a}\down{d}} 
\tensor{\CY}{\up{c}\up{d}\down{e}} 
\tensor{\V}{\down{b}\up{e}}
%
%
-224 
\tensor{\CY}{\down{a}\down{c}\down{b}} 
\tensor{\V}{\down{d}\down{e}} 
\tensor{\V}{\up{c}\up{d}\down{|}\up{e}}
\\
%
%
-96 
\tensor{\CY}{\down{a}\down{b}\down{c}} 
\tensor{\CY}{\down{d}\up{c}\down{e}} 
\tensor{\V}{\up{d}\up{e}}
%
%
-320 
\tensor{\CY}{\down{c}\down{a}\down{d}} 
\tensor{\CY}{\down{e}\down{b}\up{d}} 
\tensor{\V}{\up{c}\up{e}}
%
%
+736 
\tensor{\CY}{\down{c}\down{a}\down{d}} 
\tensor{\CY}{\up{d}\down{b}\down{e}} 
\tensor{\V}{\up{c}\up{e}}
\\
%
%
-96 
\tensor{\CY}{\down{a}\down{c}\down{d}} 
\tensor{\V}{\down{b}\up{d}\down{|}\down{e}} 
\tensor{\V}{\up{c}\up{e}}
%
%
-96 
\tensor{\CY}{\down{c}\down{a}\down{d}} 
\tensor{\V}{\down{b}\up{d}\down{|}\down{e}} 
\tensor{\V}{\up{c}\up{e}}
%
%
-192 
\tensor{\CY}{\down{c}\down{a}\down{d}} 
\tensor{\CY}{\up{c}\down{b}\down{e}} 
\tensor{\V}{\up{d}\up{e}}
\\
%
%
+16 
\tensor{\V}{\down{c}\down{d}} 
\tensor{\V}{\up{c}\down{e}} 
\tensor{\Weyl}{\down{a}\down{i}\up{d}\down{j}} 
\tensor{\Weyl}{\down{b}\up{i}\up{e}\up{j}}
%
%
-32 
\tensor{\V}{\down{c}\down{d}} 
\tensor{\V}{\up{c}\down{e}} 
\tensor{\Weyl}{\down{a}\down{i}\up{d}\down{j}} 
\tensor{\Weyl}{\down{b}\up{j}\up{e}\up{i}}
%
%
-32 
\tensor{\V}{\down{c}\down{d}} 
\tensor{\V}{\down{e}\down{i}} 
\tensor{\Weyl}{\down{a}\down{j}\down{b}\up{c}} 
\tensor{\Weyl}{\up{d}\up{e}\up{i}\up{j}}
\\
%
%
+4 
\tensor{\bg}{\down{a}\down{b}} 
\tensor{\V}{\down{c}\down{d}} 
\tensor{\V}{\down{e}\down{i}} 
\tensor{\Weyl}{\up{c}\down{j}\up{e}\down{k}} 
\tensor{\Weyl}{\up{d}\up{i}\up{j}\up{k}}
%
%
-4 
\tensor{\bg}{\down{a}\down{b}} 
\tensor{\V}{\down{c}\down{d}} 
\tensor{\V}{\down{e}\down{i}} 
\tensor{\Weyl}{\up{c}\up{e}\down{j}\down{k}} 
\tensor{\Weyl}{\up{d}\up{j}\up{i}\up{k}}
%
%
-32 
\tensor{\V}{\down{c}\down{d}} 
\tensor{\V}{\down{e}\down{i}} 
\tensor{\V}{\up{e}\up{i}} 
\tensor{\Weyl}{\down{a}\up{c}\down{b}\up{d}}
\\
%
%
+32 
\tensor{\V}{\down{c}\down{d}} 
\tensor{\V}{\up{c}\down{e}} \J 
\tensor{\Weyl}{\down{a}\up{d}\down{b}\up{e}}
%
%
-224 
\tensor{\V}{\down{c}\down{d}} 
\tensor{\V}{\down{e}\down{i}} 
\tensor{\V}{\up{c}\up{e}} 
\tensor{\Weyl}{\down{a}\up{d}\down{b}\up{i}}
%
%
+150 
\tensor{\bg}{\down{a}\down{b}} 
\tensor{\V}{\down{c}\down{d}} 
\tensor{\V}{\down{e}\down{i}} 
\tensor{\V}{\up{e}\down{j}} 
\tensor{\Weyl}{\up{c}\up{i}\up{d}\up{j}}
\\
%
%
+150 
\tensor{\bg}{\down{a}\down{b}} 
\tensor{\V}{\down{c}\down{d}} 
\tensor{\V}{\down{e}\down{i}} 
\tensor{\V}{\up{c}\down{j}} 
\tensor{\Weyl}{\up{d}\up{e}\up{i}\up{j}}
%
%
-32 
\tensor{\V}{\down{a}\down{c}} 
\tensor{\V}{\down{d}\down{e}} 
\tensor{\V}{\up{c}\down{i}} 
\tensor{\Weyl}{\down{b}\up{d}\up{e}\up{i}}
%
%
-64 
\tensor{\V}{\down{a}\down{c}} 
\tensor{\V}{\down{d}\down{e}} 
\tensor{\V}{\up{d}\down{i}} 
\tensor{\Weyl}{\down{b}\up{e}\up{c}\up{i}}
%
%
\end{array}
\]
\caption{A tensor $\SymTen_{ab}$ such that
$\cB^{8}_{ab}=\frac{1}{384}\SymTen_{(ab)}$}
\label{B8}
\end{figure}
(For typesetting convenience, the figure uses $B_{ab|cd}$ as an
alternative notation for $\nd_d\nd_c B_{ab}$, and so forth.)  To see
that $\cB_{ab}^{8}=\SymTen_{(ab)}$, we begin by constructing a
tractor formula for $\fb_{2}$ on $\ct^{2,2}[-2]$.  Let
$T\in\act^\Phi(-2)$ be an extension of any element of $\ct^{2,2}[-2]$.
By Proposition~\ref{Lab13July3g} we may assume that $\afl T=O(Q)$.
Thus by \nn{lab23June3g}, \nn{aflaR}, and \nn{abracketalX}, we first
obtain
\begin{equation}\label{L6May4a}
\begin{array}{lll}
\lefteqn{\afl \D_A T_{BCDE}=}
&&
\\
&&
-\X_A\afl^{2}T_{BCDE}
+2[\al,\nda_A]T_{BCDE}
+\aR\#\#\nda_A T_{BCDE}
\\
&&
-\frac{1}{4}\X_{A}(\aR\#\#\aR)\#\# T_{BCDE}
+ \X_{A}(\nda_{|I|}\aR)\#\# \nda^{|I|}T_{BCDE}
\\
&&
+O(Q)
.
\end{array}
\end{equation}
Here the $|\cdot|$ indicates that the enclosed index is not involved
in the hash action, and $(\aR\#\#\aR)\#\#T_{BCDE}$ denotes the double
hash of $\aR\#\#\aR$ with $T_{BCDE}$.  {}From \nn{L6May4a} together with
\nn{XRQ}, \nn{lab23June3g}, \nn{aflaR}, and \nn{L17May01c}, it follows
that
\[
\begin{array}{lll}
\lefteqn{\afl \D_A T_{BCDE}=}
&&\\
&&
-\X_A\afl^2T_{BCDE}
-2\aR_A{}^P{}_B{}^Q\D_P T_{QCDE}
-2\aR_A{}^P{}_C{}^Q\D_PT_{BQDE}
\\&&
-2\aR_A{}^P{}_D{}^Q\D_PT_{BCQE}
-2\aR_A{}^P{}_E{}^Q\D_PT_{BCDQ}
\\&&
+\frac{1}{2}\aR\#\#\D_{A}T_{BCDE}
-\frac{1}{4}\aR\#\#\X_A\aR\#\#T_{BCDE}
\\&&
-\frac{1}{4}\X_{A}(\aR\#\#\aR)\#\#T_{BCDE}
+\frac{1}{4}\X_{A}(\D_{|I|}\aR)\#\#\D^{|I|}T_{BCDE}
\\&&
-\frac{1}{8}\X_A\X_{I}(\aR\#\#\aR)\#\#\D^{|I|}T_{BCDE}
\\&&
-\frac{1}{8}\X_{A}\X^{I}(\D_{|I|}\aR)\#\#\aR\#\#T_{BCDE}
+O(Q)
.
\end{array}
\]
Since the dimension is 8, it follows from \nn{lab23June3g} that
$\X^A\D_AV=-4V+O(Q)$ for all $V\in\act^{\Phi}(-2)$.  Thus from the
definition of $\afl$, we see that
\begin{equation}\label{L8May4a}
\begin{array}{lll}
\lefteqn{
\X_A\afl^2T_{BCDE}
=}
&&\\
&&
-\al \D_A T_{BCDE}
-2\aR_A{}^P{}_B{}^Q\D_P T_{QCDE}
-2\aR_A{}^P{}_C{}^Q\D_PT_{BQDE}
\\&&
-2\aR_A{}^P{}_D{}^Q\D_PT_{BCQE}
-2\aR_A{}^P{}_E{}^Q\D_PT_{BCDQ}
\\&&
-\frac{1}{4}\aR\#\#\X_A\aR\#\#T_{BCDE}
+\frac{1}{4}\X_{A}(\aR\#\#\aR)\#\#T_{BCDE}
\\&&
+\frac{1}{4}\X_{A}(\D_{|I|}\aR)\#\#\D^{|I|}T_{BCDE}
+\frac{1}{2}\X_{A}\aR\#\#\aR\#\#T_{BCDE}
\\&&
+O(Q)
.
\end{array}
\end{equation}
We restrict \nn{L8May4a} to $\cq$ and then attach $Y^A$.
The result is that for any $T\in\ct^{2,2}[-2]$,
\begin{equation}\label{L9May4a}
\begin{array}{lll}
\lefteqn{
\fb_2 T_{BCDE}
=}
&&\\
&&
-Y^A\Box D_A T_{BCDE}
-\frac{1}{2}Y^AW_A{}^P{}_B{}^QD_P T_{QCDE}
\\&&
-\frac{1}{2}Y^AW_A{}^P{}_C{}^QD_PT_{BQDE}
-\frac{1}{2}Y^AW_A{}^P{}_D{}^QD_PT_{BCQE}
\\&&
-\frac{1}{2}Y^AW_A{}^P{}_E{}^QD_PT_{BCDQ}
-\frac{1}{64}Y^{A}W\#\#X_AW\#\#T_{BCDE}
\\&&
+\frac{1}{64}(W\#\#W)\#\#T_{BCDE}
+\frac{1}{16}(D_{|I|}W)\#\#D^{|I|}T_{BCDE}
\\&&
+\frac{1}{32}W\#\#W\#\#T_{BCDE}
.
\end{array}
\end{equation}

We next use \nn{L9May4a} to construct a tractor formula for
$\cB^{8}_{ab}$.  {}From Theorem~\ref{main2} we have
\begin{equation}\label{L10May4a}
\cB^{8}_{ce}=
{\textstyle-\frac{1}{384}}Y^BZ^C{}_c Y^D Z^E{}_e\,\fb_2W_{BCDE}.
\end{equation}
A short computation shows that
\[
W\#\#W\#\#W_{BCDE}=(W\#\#W)\#\#W_{BCDE}.
\]
Thus from \nn{L9May4a} and \nn{L10May4a} we have
%
%
\begin{equation}\label{obsdim8}
\begin{array}{lll}
\textstyle
\lefteqn{\cB^{8}_{ab}=
{\textstyle\frac{1}{24576}}Y^BZ^C{}_aY^DZ^E{}_b\left(
64 Y^A\Box D_A W_{BCDE}
\rule{0mm}{4mm}\right.}&&
\\
&&
+32Y^AW_A{}^P{}_B{}^QD_P W_{QCDE}
+32Y^AW_A{}^P{}_C{}^Q D_PW_{BQDE}
\\
&&
+32Y^AW_A{}^P{}_D{}^Q D_P W_{BCQE} +32Y^AW_A{}^P{}_E{}^Q D_P
W_{BCDQ}
\\
&&
+Y^{A}W\#\#X_AW\#\#W_{BCDE}
-3W\#\#W\#\#W_{BCDE}
\\
&&
-4(D_{|I|}W)\#\#D^{|I|}W_{BCDE}
\left.\rule{0mm}{4mm}\right)
.
\end{array}
\end{equation}
By using the same techniques as in our derivation of \nn{B6}, we see
that $\cB^{8}_{ab}=\SymTen_{(ab)}$.
%
%

%

\begin{thebibliography}{XX}
%
\bibitem{Anderson} M.T. Anderson, {\em Orbifold compactness for spaces
  of Riemannian metrics and applications}.  Math.\
  Annalen.\ \textbf{331} (2005), 739--778.  Preprint math.DG/0312111,
  http://arXiv.org.

\bibitem{Bach} R.\ Bach, {\em Zur Weylschen Relativit\"atstheorie und
    der Weylschen Erweiterung des Kr\"ummungstensorbegriffs}. Math. Z.
  \textbf{9} (1921), 110--135.
  
\bibitem{BEGo} T.N. Bailey, M.G. Eastwood, and A.R. Gover,
  \textit{Thomas's structure bundle for conformal, projective and
    related structures}.  Rocky Mountain J. Math.\ \textbf{24} (1994),
  1191--1217.

\bibitem{BEGr} T.N. Bailey, M.G. Eastwood, and C.R. Graham,
  \textit{Invariant theory for conformal and CR geometry}.  Ann.\ 
  Math.\ (2) {\bf 139} (1994), 491--552.
  
\bibitem{BC} B.D. Boe and D.H. Collingwood, \textit{A comparison
    theory for the structure of induced representations}. J.\ Alg.
  {\bf 94} (1985), 511--545.
    
\bibitem{BrSrni} T. Branson, \textit{Q-curvature and spectral
invariants}.  Supp.\ Rend.\ Circ.\ Matem.\ Palermo \textbf{75} (2005),
11--55.

\bibitem{tomsharp} T. Branson, \textit{Sharp inequalities, the functional
    determinant, and the complementary series}.  Trans.\ Amer.\ Math.\
    Soc.\ {\bf 347} (1995), 3671--3742.
  
\bibitem{BrGoPacific} T.\ Branson and A.R.\ Gover, \textit{Conformally
    invariant non-local operators}. Pacific J.\ Math.\ {\bf 201}
  (2001), 19--60.

\bibitem{BrGodeRham} T.\ Branson and A.R.\ Gover,
\textit{Conformally invariant operators, differential forms, 
cohomology and a generalisation of Q-curvature}.
Comm.\ Partial Differential Equations (to appear).
Preprint
math.DG/0309085, http://arXiv.org.

\bibitem{torsionprogress}  T.\ Branson and A.R.\ Gover, 
{\em Detour torsion}.  In progress.

\bibitem{BrGoopava} T.\ Branson and A.R.\ Gover, 
{\em Electromagnetism, metric deformations,
    ellipticity and gauge operators on conformal 4-manifolds}.
Diff.\ Geom.\ and its Applications {\bf 17} (2002), 229--249. 

\bibitem{tbbo91} T. Branson and B. {\O}rsted, {\em 
Explicit functional determinants 
in four dimensions}.  Proc.\ Amer.\ Math.\ Soc.\ {\bf 113} (1991),
669--682.

\bibitem{CapGoamb} A. \v Cap and A.R. Gover, {\em Standard tractors
    and the conformal ambient metric construction}. Ann. Global Anal.
  Geom. \textbf{24} (2003), no. 3, 231--259.  Preprint
  math.DG/0207016, http://arXiv.org.
  
\bibitem{CapGoluminy} A. \v Cap and A.R. Gover, \textit{Tractor
    bundles for irreducible parabolic geometries}.  Global analysis
  and harmonic analysis (Marseille-Luminy, 1999), 129--154, S\'emin.
  Congr., \textbf{4}, Soc. Math. France, Paris, 2000.  Preprint
  Publications/SeminairesCongres, http://smf.emath.fr.

\bibitem{CapGotrans} A. \v Cap and A.R. Gover, \textit{Tractor calculi for
parabolic geometries}. Trans. Amer. Math. Soc. \textbf{354} (2002),
1511--1548.  Preprint 792, http://www.esi.ac.at.

\bibitem{Esrni} M.G. Eastwood, {\em Notes on conformal differential
    geometry}.  The Proceedings of the 15th Winter School ``Geometry
  and Physics'' (Srni, 1995). Rend. Circ. Mat. Palermo (2) Suppl. No.
  43 (1996), 57--76.
  
\bibitem{EastSlo} M.G. Eastwood and J. Slov\'ak, Semiholonomic Verma
  modules. {\em J.\ Algebra} {\bf 197} (1997), 424--448.

\bibitem{F} C. Fefferman, \textit{Parabolic invariant theory in
    complex analysis}. Adv.\ in Math.\ {\bf 31} (1979), 131--262.
  
\bibitem{FGast} C. Fefferman and C.R. Graham, \textit{Conformal
    invariants}.  The mathematical heritage of \'Elie Cartan (Lyon,
  1984). Ast\'erisque 1985, Numero Hors S\'{e}rie, 95--116.

\bibitem{FeffGrinprogess} C.\ Fefferman and C.R.\ Graham, {in progress}.
  
\bibitem{GG} J. Gasqui and H. GoldSchmidt, {\it D\'eformations
    infinit\'esimales des structures conformes plates (French)
    [Infinitesimal deformations of flat conformal structures]}.
  Progress in Mathematics, 52. Birkh\"auser Boston, Inc.,
  Boston, MA, 1984.

\bibitem{Gau} P.\ Gauduchon, {\em Connexion canonique et structures de
    Weyl en g\'{e}om\'{e}trie conforme}. Preprint, CNRS UA766, Ecole
    Polytechnique (1990).

\bibitem{gosrni} A.R. Gover, \textit{Aspects of parabolic invariant theory}.
The 18th Winter School ``Geometry and Physics'' (Srni, 1998). Rend.
   Circ. Mat. Palermo (2) Suppl. No. 59 (1999), 25--47. 

\bibitem{goadv} A.R. Gover, \textit{Invariant theory and calculus for
    conformal geometries}. Adv.\ Math.\ {\bf 163} (2001), 206--257.

\bibitem{GoH} A.R.\ Gover and K.\ Hirachi, {\em Conformally
  invariant powers of the Laplacian -- A complete non-existence theorem}.
J. Amer. Math. Soc.\ {\bf 17}  (2004), 389--405. 
Preprint math.DG/0304082, http://arXiv.org.

\bibitem{GoNur} A.R.\ Gover and P.\ Nurowski,
{\em Obstructions to conformally Einstein metrics in $n$ dimensions}.
J. Geom.\ Phys.\ (to appear).
  Preprint, math.DG/0405304, http://arXiv.org.
  
\bibitem{GoPet} A.R. Gover and L.J. Peterson, \textit{Conformally
    invariant powers of the Laplacian, Q-curvature, and tractor
    calculus}. Comm. Math. Phys.\ \textbf{235} (2003), no. 2, 339--378.
  Preprint math-ph/0201030, http://arXiv.org.

\bibitem{Grpriv} C.R. Graham, {\em private communication}.

\bibitem{GrH} C.R. Graham and K. Hirachi, {\em The ambient obstruction
tensor and $Q$-curvature}.  AdS-CFT Correspondence: Einstein Metrics
and their Conformal Boundaries. IRMA Lect.\ Math.\ Theor.\ Phys.\ {\bf
8} (to appear).  Preprint math.DG/0405068, http://arXiv.org.

\bibitem{GJMS} C.R. Graham, R. Jenne, L.J. Mason, G.A. Sparling,
  \textit{Conformally invariant powers of the Laplacian, I:
    Existence}. J. London Math. Soc.\ \textbf{46} (1992), 557--565.

\bibitem{Lee} J. Lee, \textit{Ricci} software package,
http://www.math.washington.edu/$\tilde{\ }$lee.

\bibitem{TianViac} G.\ Tian and J.\ Viaclovsky, {\em Bach-flat
    asymptotically locally Euclidean metrics}.  Invent.\ Math.\
    \textbf{160} (2005), 357--415.  Preprint math.DG/0310302,
    http://arXiv.org.

\bibitem{T} T.Y. Thomas, {\em On conformal geometry}.
  Proc.\ Natl.\ Acad.\ Sci.\ USA
  {\bf 12} (1926), 352--359.

\end{thebibliography}
\end{document}